\documentclass[12pt]{article}
\usepackage{amssymb}
\usepackage{latexsym}
\usepackage{exscale}

\def \beq{\begin{equation}}
\def \eeq{\end{equation}}
\def \lab{\label}
\renewcommand{\rq}[1]{(\ref{#1})}
\newtheorem{lemma}{Lemma}
\newtheorem{prop}{Proposition}
\newtheorem{thm}{Theorem}
\newtheorem{cor}{Corollary}

\newcommand{\bR}{{ \mathbb R  }}
\newcommand{\bC}{\Bbb C}
\newcommand{\bZ}{\Bbb Z}
\newcommand{\bN}{\Bbb N}
\newcommand{\bK}{\Bbb K}

\newcommand{\E}{\mathcal{E}}

\newcommand{\la}{\mbox{$\lambda$}}

\newcommand{\ra}{\mbox{$\mapsto $}}

\newcommand{\ep}{\epsilon}

\newcommand{\N}{{\cal N}}
\newcommand{\al}{\alpha }

\newcommand{\ga}{\gamma }
\newcommand{\f}{\varphi }

\newcommand{\La}{\Lambda }

\newcommand{\mC}{{\mathcal C}}
\begin{document}
\title{Exact controllability for string  with attached masses}

{\bf accepted by SIAM Journal of Optimization and Control}

 \maketitle

\noindent{\bf Sergei Avdonin}, Department of Mathematics and Statistics, University of Alaska at Fairbanks, Fairbanks,
AK 99775, U.S.A.; s.avdonin@alaska.edu; (907)-474-5023

\

\noindent{\bf Julian Edward}, Department of Mathematics and Statistics, Florida International University, Miami,
FL 33199, U.S.A.; edwardj@fiu.edu; (305)-348-3050

\

{\bf Abstract:} We consider the problem of boundary control for a vibrating string
with $N$  interior point masses. We assume the control is at the left end, and the string is fixed at the right end.
Singularities in waves are ``smoothed" out to one order as they cross
a point mass. We characterize the reachable set
 for  a $L^2$ control. The control problem is reduced to a moment problem, which is then solved
 using the theory of exponential divided differences in tandem with  unique shape and velocity controllability results.

\begin{section}{Introduction}

There has been much  interest in so called ``hybrid systems" in which the dynamics of elastic systems
and possibly rigid structures are related through some form of coupling. The study of
controllability and stabilization of such structures has made in a number of works, see
\cite{LLS} and  \cite{DZ} and references therein, also \cite{HZho}.
Networks of strings with attached masses have also been studied by
 many  of authors in the context of inverse problems,  see for instance \cite{Eck} and references therein.
 The controllability of a string with a single attached mass was considered in
 \cite{HZ}, \cite{C}, \cite{CZ}. The controllability of a series of Euler-Bernoulli beams with interior attached masses
 was considered in \cite{MR}.

 We consider the wave equation on the interval $ [0,\ell]$  with masses $M_j >0$ attached at the  points $a_j, \;  j=1,\ldots,N, $
where $ 0=a_0< a_1 < ... < a_N<a_{N+1}=\ell.$ We will assume the string and masses are at rest and  at equilibrium  until time $t=0$, when a control is applied
at the left end of the string with the right end fixed.
In what follows, $u_j(x,t)$ will denote the vertical displacement of the string in the interval
$(a_j,a_{j+1})$, and $h_j(t)$ will denote the vertical displacement of the mass at $x=a_j$.
On each interval $(a_j,a_{j+1})$, let  $q_j (x)$ represent   some potential.
Small amplitude vibration of this system
is modeled by
\begin{eqnarray}
 \frac{\partial ^2u_j}{\partial t^2}- \frac{\partial ^2u_j}{\partial x^2}+q_{j}(x)u_j& =&0, \ t\in \bR,\ x\in (a_j,a_{j+1}),\ j=0,\ldots,N,\label{A_N11}\\
u_j(x,t)& =& 0,\ t\leq 0,\ j=0,\ldots,N,\\
h_j(x,t)& =& 0,\ t\leq 0,\ j=0,\ldots,N,\\
u_{j-1}(a_j^-,t)& =& h_j(t)=u_j(a_j^+,t), \ j=1,\ldots,N, \label{contu}\\
{M_j} h_j''(t)& = & \frac{\partial u_j}{\partial x}(a_j^+,t)-\frac{\partial u_{j-1}}{\partial x}(a_j^-,t),\ j=1,\ldots,N,\label{Newton}\\
  u_{N}(\ell^- ,t)& =& 0,
 \end{eqnarray}
 with the control,
\begin{equation}
u_0(0^+,t) = f(t). \label{control9}
\end{equation}
Here $u(a_j^+,t):=\lim_{\ep \to 0^+}u(a_j+\ep ,t)$ for fixed $t$, and $u(a_j^-,t)$ is defined similarly.

  We assume for $j=0,1,2,$ that $q_j$ extends to
  	 $ C[a_j,a_{j+1}]$, while for $j> 2$,  $q_j$ extends
to a function in
 $ C^{j-2}[a_j,a_{j+1}]$. Our methods still apply
  if $q_j\in
  H^{\max (0,j-2)}(a_j,a_{j+1})$, but the presentation more cumbersome.
Here and it what follows $H^j(a_j,a_{j+1})$ refers to the standard Sobolev space, with $H^0 =L^2$.
Define
$\theta ^{-1}(0,a_1):=\{ u\in H^1(0,a_1): u(0)=0\} '$.
 We define
$$\tilde{W}_0=\oplus_{j=0}^NH^{j}(a_j,a_{j+1}) \mbox{ and } \tilde{W}_{-1}=\theta^{-1}(0,a_1)\oplus \big ( \oplus_{j=1}^{N}H^{j-1}(a_j,a_{j+1}) \big ).$$
One of the most important features of System \rq{A_N11}-\rq{control9} is that
 the attached
masses will mollify transmitted waves, so the system  is well posed in asymmetric spaces.
This is reflected in the following:
\begin{prop}\label{Nreg0} For any $T>0$, let $f\in L^2(0,T)$. There exists a unique solution
$$(u_0,h_1,u_1,h_2, ...,h_N, u_N)$$
to System \rq{A_N11}-\rq{control9}.
 For $i=0,1$,
\begin{equation}
(u_0,u_1, ..., u_N)\in C^i([0,T];\tilde{W}_{-i}).\label{spacescaled}
\end{equation}
Furthermore, $h_j(t) \in H^j_{loc}$ for each $j$.
\end{prop}
Since $f\in L^2$, the vector $(u_0,h_1,...,h_N,u_N)$ will actually be a weak solution to the system.
To state our controllability results, it is convenient to disregard for the moment the states
 of the masses.  We will discuss those  after the statement of Theorem \ref{mainthm}.
In what follows, we will refer to the vector $(u_0,u_1,..., u_N)$, simply as $u^f(x,t)$.
We say the pair of functions $(\phi_0 (x),\phi_1(x))$ is in the ``reachable set at time $T$'' if there exists $f\in L^2(0,T)$ such that $(u^f(x,T),u_t(x,T))=(\phi_0(x),\phi_1(x))$ for $x\neq a_j$. We wish to  characterize the reachable sets. For each $j$, the masses impose on $(u^f(x,T),u^f_t(x,T))$ a set of equations that must hold at $x=a_j$, provided
$u^f$ and $u^f_t$ are sufficiently regular. One example of this is $u^f(a_j^-,T)=u^f(a_j^+,T)$, which by \rq{contu} and Proposition \ref{Nreg0} must hold for all $j\geq 2$.
In addition, the boundary condition at $x=\ell$ imposes further conditions.
The collection of such equations satisfied by $u^f(x,T)$ will be denoted   ${\mathcal C}^{0}_*$, while
 the collections of equations for $u^f_t(x,T)$ will be denoted   ${\mathcal C}^{-1}_*$. These spaces will be carefully
  described in Section \ref{solutionrep}.

 We now define a Hilbert space ${W}_i$,  for integers  $i=-1,0$, by
$${W}_i=\{ \phi(x)\in  \tilde{W}_i: \ \phi \in {\mathcal C}^{i}_*\}.$$

We  have the following
\begin{thm}\label{reachable}
 For any $ T > 0 , $
 $$ \{ (u^f(\cdot,T),u^f_t(\cdot ,T)) : f \in L^2(0,T) \} \subset W_0\times W_{-1}. $$
\end{thm}

We can now state our main result.
\begin{thm} \label{mainthm}  Let $N\geq 1$, and let $T > 2\ell$.
Then for any $(y_0,y_1)\in W_0\times W_{-1}$,
 there exists a control $f \in L^2(0,T)$ such that the solution $u^f$ to \rq{A_N11}-\rq{control9} satisfies
 \begin{equation}
  u^f(x,T)=y_0(x),\, u_t^f(x,T)=y_1(x), \ x \in (0,\ell).\label{fullcontrol}
 \end{equation}
Furthermore,
\begin{equation}
\|  f\|^2_{L^2(0,T)}\asymp  \| y_0\|^2_{W_0}+\| y_1\|^2_{W_{-1}}.\label{energyequiv}
\end{equation}
\end{thm}
Here and below,
$||f||_A\asymp ||g||_B$ means there exist positive constants $C_1,C_2$ such that
$C_1||f||_A \leq ||g||_B\leq C_2||f||_B$  for all  $f \in A,\, g \in B$. In what follows, we will refer
to the conclusions of Theorem \ref{mainthm} as ``full controllability".

Theorem \ref{mainthm} makes assertions about the terminal positions and velocities of the string segments, but not directly about the terminal positions and velocities of the masses, which we discuss now.
By Proposition \ref{Nreg0}, \rq{contu}, and \rq{Newton} we have
\begin{equation}
h_j(T)=u^f(a_j^+,T), \ j\geq 1, \ \mbox{ and } \ h_j'(T)=u_t^f(a_j^+,T), \ j\geq 2.\label{m-s}
\end{equation}
Thus the $\{ h_j(T),h_j'(T)\}$ are determined by the terminal string positions and velocities except for $h_1'(T)$. In this paper we will show $h_1'(t)$ is $L^2$ but
not necessarily continuous for $f\in L^2$, so $h_1'(T)$ cannot generally be  prescribed as terminal data.

 We now state some results that are used in the proof of Theorem \ref{mainthm}, and that
we believe are of independent interest.
 Let
 $$W^T_i:=\{ \phi \in W_i: \phi(x)=0,\ \forall x\geq T \} ,\ \, T \leq \ell, \ i=-1,0 \ \mbox{ and } \
 W^T_i=W_i \ \mbox{ for } \ T>\ell.$$

\begin{thm} \label{shapecont}

Let $T>0$.

\noindent A)   For any  $\phi \in W_0^T$, there exists  $f\in L^2(0,T)$
such that $u^f(x,T)=\phi(x)$ and
\begin{equation}
\| u^f(\cdot,T)\|_{W_0} \asymp \| f\|_{L^2(0,T)}.\label{sheq}
\end{equation}
If $T\leq \ell$, this $f$ is unique.

\noindent B) For any  $\phi \in W_{-1}^T$, there exists a $f\in L^2(0,T)$
such that $u^f_t(x,T)=\phi(x)$ and
\begin{equation}
\| u^f_t(\cdot,T)\|_{W_{-1}} \asymp \| f\|_{L^2(0,T)}.\label{veleq}
\end{equation}
If $T\leq \ell$, this $f$ is unique.

\end{thm}
That an arbitrary position (resp. velocity) within the appropriate function space can be attained in time $T $ will be referred to as ``shape controllability in time $T\,$" (resp. ``velocity controllability in time $T\, $").  Shape and velocity controllability defined as above are  key ingredients
of the boundary control method in inverse theory (see, e.g. \cite{ABI,ALP}). The main tool used to prove Theorem \ref{shapecont} is a representation
of $u^f$ for $x\in (0,a_1)$ which can be thought of as a perturbation of the d'Alembert  solution. This representation has been used
previously in many works on systems without masses, see for instance \cite{AK}. This  representation, together with a careful analysis of the relation
between the control $f$ and the positions of the masses, are used to reduce the shape and velocity control problems to  Volterra integral equations of the second kind. This material will be proven in Sections 2 and 3.

We now give some ideas of the proof of Theorem \ref{mainthm}. The spectral problem associated to System \rq{A_N11},\rq{control9} is
\begin{eqnarray}
-\phi''(x)+q(x)\phi (x) & = & (\la)^2\phi (x),\ x\in (0,\ell )\setminus \cup_{j=1}^Na_j,\nonumber \\
\phi (0)& = & \phi (\ell)  =  0,\nonumber\\
\phi (a_j^-)& =&\phi (a_j^+),\nonumber\\
\phi '(a_j^+) & = & \phi '(a^-_j)-M_j\la^2\phi (a_j^+),\ j=1,\ldots , N. \label{SL}
\end{eqnarray}
Let the  lengths of the subintervals $[a_{j},a_{j+1}]$ be given by
 $ \{ \ell_j\}_0^N.$ Let $\{ (\la_n)^2: \ n \in \bN\}$ be the set of eigenvalues for the system \rq{SL}, listed in increasing order.
Taking (possibly complex) square roots, we then define the associated eigenfrequencies
 $\Lambda :=\{ \la_n: \ n =\pm 1, \pm2, ... \}$.
 In Section 4, we prove
   \begin{thm}\label{spectrumq}
Let $\Lambda '$ be any subset of $\Lambda$ obtained by deleting $2N$ elements. Then $\Lambda '$ can be reparametrized as
$$\Lambda ' = \cup_{j=0}^N\{ \la^{(j)}_m\}_{m\in \bK},$$
where for each $j$,
 \begin{equation}
 |\la_m^{(j)}-\frac{\pi m}{\ell_j}|=O(|m|^{-1}).\label{babyasy2}
 \end{equation}
 \end{thm}
We use the $\Lambda $
 to construct an associated family of exponential  divided differences (E.D.D.) using ideas developed in \cite{AI1}. Fix $T>\ell$.
Applying a result found in \cite{AM1} to  Theorem \ref{spectrumq},
we show that this family
 forms a Riesz sequence in $L^2(0,2T)$.
 We can then solve the moment problem associated to \rq{fullcontrol}, but only in terms of a space of Fourier coefficients whose characterization is not obvious.
The novelty in our argument is in what follows next in Section 5.
Corresponding to the E.D.D.  are divided differences of sines and cosines, which by \cite{ABI} form
 Riesz sequences
	 on $L^2(0,T )$.
This, together with the estimates
\rq{sheq} and \rq{veleq}, allow us to construct Riesz bases of $W_0$ and $W_{-1}$ from the eigenfunctions associated to \rq{SL}.
After this, we complete the proof by proving \rq{energyequiv}.

This paper was inspired by the paper of Hansen and Zuazua, \cite{HZ}.
In that work (also see \cite{CZ}), the authors consider only the case $N=1$. Their main result  is that for $T>2\ell$,
the reachable set is $W_0\times W_{-1}$.
 Their method of proof
involves using the theory of characteristics for the constant coefficient wave equation to prove an observability estimate.
They then suggest how their results can be extended to the case of variable coefficients. It is not clear that the methods
in their paper can easily be extended to our setting with many masses, and they state that this case is beyond the scope of their paper.
Our approach is quite different from \cite{HZ}. Ours  combines dynamical (Sections 2, 3) and spectral methods
	(Sections 4, 5) that allow us to split  the problem to two subproblems, each of them is interesting and nontrivial. The first one, shape and velocity controllability, is solved by dynamical methods. The main achievement of our spectral investigations
	is the successful application of the method of moments to systems with nonsymmetric reachable spaces. In addition, our construction of Riesz basis
for asymmetric spaces is applicable in other settings, some of which we mention below.
%\footnote{HZ almost achieve this in the sense that
%they present a Riesz basis for $W_0$ in the special case $a=\ell/2, N=1, q=0$; see Lemma 5.6 there. However, they don't note the significance of this. %Is it worth commenting on this?
%{\bf S: I think that it would too complicated and unclear idea for a standard reader; it is enough what we say in the next paragraph}}, which is based %on the modern theory of Riesz bases of EDD \cite{AI1,AM1}.

An interesting question is whether  or not the time for full control can be improved to $T=2\ell$. This topic will be discussed in a
forthcoming paper,\cite{AE2}, where we show that full control can be achieved in time $T=2\ell$ for $N=1$, but not for $N>1$. In that paper,
we also consider the case of Neumann or mixed control at one end with the other end fixed, and with variable density and tension. Unlike the case in this paper, in that case
one can have a uniform gap in the spectral frequencies, but the method of using E.D.D. is still useful in that case. In that paper, we also
construct a scale of asymmetric Sobolev-like spaces starting from $W_0$ and $W_{-1}$, and we study their properties.

The results obtained in this paper and \cite{AE2} can be extended to networks of strings with loaded masses,
see \cite{AAE} for tree-like networks. Our results on Riesz bases of non-symmetric spaces also can be applied to the controllability
 of beam and dynamical Schr\"{o}dinger equations. This will be a topic of a forthcoming paper.

% Hansen and Zuazua also study controllability
%using spectral analysis in the special case where $q=0$ and $a_1=\ell/2$, and they are able to given a rather precise
%spectral characterization of the reachable set. However, as we observe in this paper, it might not be easy to extend this
%analysis to a more general case. In particular, in the case $q=0$ and $a_1=\ell/2$ and $N=1$, the associated unit norm eigenfunctions
%have the asymptotics $|\phi_j'(0)|\asymp j$, whereas we shall show that this is not generally the case.

This paper is organized as follows. In Section 2 first we prove the  regularity results, and give a representation of the solution $u$ that effectively models the propagation, transmission, and
reflection of waves. We then discuss the various compatibility conditions that arise. In Section 3, we prove the shape and velocity controllability results. In Section 4, we study the spectral theory associated to
System \rq{A_N11}-\rq{control9} to obtain a Riesz sequence result for an associated family of  E.D.D. Then in Section 5, we prove
the main result.

During the  submission of this manuscript, a manuscript by Ben Amara and Beldi  came to our attention, \cite{AB1}. There, the authors consider a vibrating string with one attached mass, variable density, tension and $q_j$, and Dirichlet control applied at one end with the other  end fixed, and  they prove exact controllability. Their proof uses a precise analysis of the eigenvalue and eigenfunction asymptotics, and a version of Ingham's inequality that applies in weakened gap conditions \cite{BKL}. It is not clear that their methods would apply in the N mass case.

\vskip2mm
{\bf Acknowledgements.} The authors gratefully thank Scott Hansen for useful conversations. The research of Sergei Avdonin was supported in part by the National Science Foundation,
grant DMS 1411564 and by the Ministry of Education and Science of Republic of Kazakhstan, grant no. 4290/GF4. The authors would also like to thank
the referees for their comments that improved our exposition.

\end{section}
\begin{section}{Existence, uniqueness, and regularity of solutions}
In this section, we study  properties of solutions of  the system \rq{A_N11}-\rq{control9}. The case $N=1$ was studied in \cite{HZ},
but for larger $N$ a more detailed analysis is required for various reasons, for instance  because the reachable set must satisfy
many compatibility conditions at the masses.
 Our method
of solving for $u$ is  different from \cite{HZ}, and uses a perturbed version of d'Alembert's solution to the wave equation commonly used in the study of the inverse problems, see for instance \cite{AK}. This representation, together with a careful analysis of the relationship between $f(t)$ and
$h_1(t)$ (see Lemma \ref{existenceh}), are the key ingredients in the proof of Theorem \ref{shapecont}. In Section 2.1, we list some notation used in this paper. Then in Section 2.2, we give a detailed study of the case $N=1$.  The arguments from Section 2.2 are then extended to general $N$ in   Section 2.3, where   we also discuss the certain equations that must be satisfied by functions in
 the reachable set and finish the proof of Theorem 1.

\begin{subsection}{Notational preliminaries}
For functions of only one variable,
 $f^{(j)}$ denotes the $j$th derivative, but for $j=1,$ resp. 2 we will without confusion use the notation
$f'$, resp. $f''$. For partial derivatives of $u=u(x,t)$, it  will be convenient to use either
 $u_{xx}$ or $\frac{\partial^2u}{\partial x^2}$.

Fix $T>0$.
Define $H^{j}(a,b)$ to be the set of functions in $L^2[a,b]$ whose
weak derivatives up to order $j$ are in $L^2[a,b]$. The corresponding norms will be denoted $||*||_{H^j(a,b)}$.
 We set the following notation
 $$ C^j_*= \{ f\in C^j(-\infty ,T]: f(t)=0 \mbox{ if } t\leq 0\},$$
 $$ L^2_*= \{ f\in L^2(-\infty ,T): f(t)=0 \mbox{ if } t\leq 0\},$$
$$ H^{n}_*= \{ f\in H^{n}(-\infty ,T): f(t)=0 \mbox{ if } t\leq 0\}.$$
Note that  $f\in H^{n}_*$ implies $f^{(j)}(0)=0$ for $j=0,\ldots ,n-1$.
We also define
$$H^{-1}_*= \{ f\in H^{-1}(-\infty ,T): f|_{(-\infty ,0)}=0 \mbox{ as a distribution. } \}.$$
Finally, when convenient we will denote the solution of System \rq{A_N11}-\rq{control9} by $u$ rather than $u^f$.

\end{subsection}
\begin{subsection}{Case of single mass}
In this subsection we give a representation of the solution to  System
\rq{A_N11}-\rq{control9} in the case of a single mass located at $x=a$, with $a \in (0,\ell)$.
This local representation
can readily be extended to the case $N>1$. Thus consider
\begin{eqnarray}
u_{tt}-u_{xx}+q(x)u& =&0, \ t\in (0,T),\ x\in (0,\ell)\setminus \{ a\},\label{pdn1}\\
u(x,t)& =& 0,\ t\leq 0,\label{inn1}\\
u(a^-,t)& =& u(a^+,t)=h(t-a), \label{bead1}\\
M h''(t-a)& = & u_x(a^+,t)-u_x(a^-,t),\label{bead2}\\
u(0,t)& =& f(t),\label{bcn1}\\
  u(\ell ,t)& =& 0.\label{SA}
\end{eqnarray}
We recall the
following  result proven in \cite{HZ}.
\begin{thm}\label{eur}
For $f\in L^2(0,T)$, there exists a unique weak solution $(u,h)$ to System \rq{pdn1}-\rq{SA}.
Denote by $u_j$ the restriction of $u(x,t)$ to $(a_j,a_{j+1})$, with $ j=0,1;\,a_1=a$. Then
$$
(u_0,u_1)\in C^i([0,T];\tilde{W}_{-i}).
$$
\end{thm}
This result is proven in \cite{HZ}, but the representation of  $u$ in this paper is  different.

Fix $b\in \mathbb{R}$.
We begin with the representation of the solution $u$. As a preliminary step,
let $D =\{ (x,t)|\ 0<x<t<\infty \}$.
 Consider the Goursat problem:
\begin{equation}
k_{tt}(x,t)-k_{xx}(x,t)+q(x+b)k(x,t)=0,\ (x,t)\in D,\label{g1}
\end{equation}
\begin{equation}
k(0,t)=0,\ k(x,x)=-\frac{1}{2}\int_{0}^{x} q(\eta +b )d\eta .
\label{g2}
\end{equation}

 \begin{prop} \label{goursat}
 \

a)  Fix $b\in \bR$. Let ${q}\in C_{loc}[b,\infty )$.
 The system (\ref{g1},~\ref{g2}) has a unique generalized
  solution, denoted $k(b^+; x,t)$, such that $k(b^+;\cdot,\cdot)\in C^{1}(\overline {D} )$ and the boundary conditions hold in a classical sense.

b) Fix $b\in \bR$ and $n\in \bN$. Let ${q}\in C^{n}_{loc}[b,\infty )$.
 The system (\ref{g1},~\ref{g2}) has a unique
  solution, denoted $k(b^+; x,t)$, such that $k(b^+;\cdot ,\cdot)\in C^{n+1}(\overline {D} )$ and \rq{g1} and \rq{g2} hold in a classical sense.

\end{prop}
For a proof of this, the reader is referred to  (\cite{AM}, Proposition 6).

We now solve the system
\begin{equation}
u_{tt}-u_{xx}+q(x)u=0, \ t\in (0,\infty ),\ x\in (b,\infty ) ,\label{pdeR}
\end{equation}
\begin{equation}
 u(b,t)=f(t),\ t>0,\label{bcR}
 \end{equation}
\begin{equation}
 u(x,t)=0,\ x>b,\ t\leq 0 .\label{icR}
 \end{equation}

The following holds by direct calculation:
\begin{prop}\label{wavesol}
Let $k$ be as in Proposition ~\ref{goursat}.

 a) Suppose $f\in C_*^2$. Then the
problem \rq{pdeR},\rq{bcR},\rq{icR} has unique solution
$u^f(b^+;x,t)$, with
\begin{equation}
u^f(b^+;x,t)=
f(t-x+b)+\int_{s=x-b}^tk(b^+;x-b,s)f(t-s)ds; \label{wf1}
\end{equation}
$u^f\in H^{2}((b,\infty)\times
(0,T))$, \rq{pdeR} is satisfied almost everywhere, and the
boundary and initial conditions are satisfied in a classical sense.

b) For $f\in L^2_*$, the function $u^f(b^+;x,t)$ defined above
gives a  solution to
\rq{pdeR} in the distribution sense,  \rq{bcR} holds for almost all $t$, and \rq{icR}  holds for all $x$.
Furthermore,   $u^f\in
C([0,T];L^2(b,\infty ))$.

c)  Let integer $j$ satisfies $j\geq 2$. Suppose $f\in H^j_*$ and ${q}\in C^{j-2}$.
Then for  fixed $t$, $x\mapsto u^f(b^+;x,t)$ is in $H^j_{loc}$, and for  fixed $x$, $t\mapsto u^f(b^+;x,t)$ is in $H^j_{loc}$.

\end{prop}
{\bf Remark on notation:} the superscript on  $b^+$ in $u^f(b^+;x,t)$
  serves to indicate that the associated wave will propagate to the right; and for
 the same reason we use the $+$ superscript on the kernel  $k(b^+; x,t)$;
similarly below $k(b^-;x,t)$ will be the kernel associated to waves propagating to the
left.

\noindent{\bf Remark:} Since $f(t)=0$ for $t<0$, \rq{wf1} can be rewritten as
$$
u^f(b^+;x,t)=\left \{
\begin{array}{cc}
f(t-x+b)+\int_{s=x-b}^tk(b^+;x-b,s)f(t-s)ds,&  x-b<t,\\
0,& x-b\geq t. \end{array} \right . $$

\

\

Setting $u(x,t)=u^f(0^+;x,t)$ as in \rq{wf1}, $u$ is the solution for System \rq{pdn1}-\rq{SA}
  $t<a$.
To consider $t\geq a$, we must study the interaction between the wave $u^f$ and the mass.
In preparation for this, note that
similar to Proposition \ref{wavesol},  the system
$$
u_{tt}-u_{xx}+q(x)u=0, \ t\in (-\infty ,\infty ),
\ x\in (-\infty ,a ) ,
$$
$$
 u(a,t)=\tilde{f}(t-a),\ t>0,\ \tilde{f}\in L^2_*,
$$
$$
 u(x,t)=0,\ x<a,\ t\leq 0.
$$
is solved by
\begin{equation}
u^{\tilde{f}}(a^-;x,t):=\tilde{f}(t+x-2a)+\int_{s=a-x}^{t-a}k(a^-;a-x,s)\tilde{f}(t-a-s)ds,
\label{reflectg}
\end{equation}
where $k(a^-;x,s)$ is the obvious analogue to $k(b^+;x,s)$.
We note
in this case
$k(a^-;0,s)=0.$

\

We now begin to solve for $u$ when $t>a$, which means the wave will be interacting with the mass.
  Define $h(t-a)=u(a,t)$; thus clearly
$h(s)=0$ for $s\leq 0$. Then for $t\in [a, \min  (2a,\ell ))$ we
have for $x>a$ by Proposition~\ref{wavesol}
\begin{equation}
u(x,t)=h(t-x)+\int_{s=x-a}^{t-a}k(a^+;x-a,s)h(t-a-s)ds;\label{bigx}
\end{equation}
this corresponds to the wave transmitted across the mass before it has reached the end at $x=\ell$.
Define $g(t)$ by
$$g(t-a)=h(t-a)-f(t-a)-\int_{s=a}^{t}k(0^+;a,s)f(t-s)ds.$$
Then for $x<a$, we have
\begin{eqnarray}
u(x,t)& =& f(t-x)+\int_{s=x}^{t}k(0^+;x,s)f(t-s)ds\nonumber \\
& & +g(t+x-2a)+\int_{s=a-x}^{t-a}k(a^-;a-x,s)g(t-a-s)ds \label{fgh} \\
&  =& f(t-x)+\int_{s=x}^{t}k(0^+;x,s)f(t-s)ds\nonumber \\
& &+ h(t+x-2a)-f(t+x-2a)-\int_{s=a}^{t+x-a}k(0^+;a,s)f(t+x-a-s)ds\nonumber \\
&
&+\int_{s=a-x}^{t-a}k(a^-;a-x,s)\left ( h(t-s-a)-f(t-s-a)-\int_{r=a}^{t-s}k(0^+;a,r)f(t-s-r)dr\right ) ds.\nonumber\\
& &\label{smallx}
\end{eqnarray}
Thus the sum of terms involving $g$  in \rq{fgh} correspond to the wave reflected off the mass. Note that by the definition of $g$,
\rq{g2} and \rq{fgh}, we have that the continuity condition, \rq{bead1}, is satisfied. The
condition \rq{bead2} implies by \rq{bigx},\rq{smallx}:
\begin{eqnarray}
M h''(t-a)&=&-2h'(t-a)+\int_{s=0}^{t-a}\frac{\partial k}{\partial
x}(a^+;0,s)h(t-a-s)ds\nonumber\\
& +& 2f'(t-a)+2k(0^+;a,a)
f(t-a)-\int_{s=a}^{t}\frac{\partial k}{\partial x}(0^+;a,s)
f(t-s)ds \nonumber\\
& +& \int_{s=a}^t\frac{\partial k}{\partial s}(0^+,a,s)f(t-s)ds\nonumber \\
& +&\int_{s=0}^{t-a}\frac{\partial k} {\partial x}(a^-;0,s) \left (
h(t-a-s)-f(t-a-s)-\int_{r=a}^{t-s}k(0^+,a,r)f(t-s-r)dr \right )ds.\nonumber \\
 & & \label{h1}
\end{eqnarray}
We now discuss the existence and regularity of the function $h$.
\begin{lemma}\label{existenceh}

Let  $T>0$.

A) Given $f\in C_*^2$ , there exists a unique $h\in
C^3_*$ solving \rq{h1} for all $t\leq T$.

B) Define the mapping $S$ by
$$
(Sf)(t)=h(t).
$$
Then $S$ is well defined, and extends to a bounded and boundedly invertible linear mapping $L^2_*\rightarrow
H^{1}_*$.

C)  $S$ extends to a bounded and boundedly invertible linear mapping $H^j_*\ra H^{j+1}_*$ for any positive
integer $j$.

\end{lemma}

Proof:  We begin with part A, so $f\in C^2(-\infty, T)$ with $f(t)=0$ for $t<0$.

We now show that a solution to \rq{h1} exists and is in $C^3$.
We rewrite \rq{h1} as
\begin{equation}
 M h''(t-a)+2h'(t-a)={\psi}(t-a )+
2f'(t-a)+\phi^2(t-a)+\phi^3 (t-a),\label{Mh}
\end{equation}
 where
\begin{eqnarray}
\phi^2(t-a)&=& 2k(0^+;a,a)f(t-a),\nonumber\\
\phi^3 (t-a) & =&
 -\int_{s=a}^{t}\frac{\partial k}{\partial x}(0^+;a,s)
f(t-s)ds
+\int_{s=a}^t\frac{\partial k}{\partial s}(0^+;a,s)f(t-s)ds\nonumber \\
&- & \int_{s=0}^{t-a}\frac{\partial k} {\partial x}(a^-;0,s) \left (
f(t-a-s)+\int_{r=a}^{t-s}k(0^+,a,r)f(t-s-r)dr \right )ds,\nonumber \\
& & \label{phi3a}
\end{eqnarray}
 and
$$
\psi (t-a)=\int_{s=0}^{t-a}\left ( \frac{\partial k}{\partial
x}(a^-;0,s)+\frac{\partial k} {\partial x}(a^+;0,s)\right )
h(t-a-s)ds.
$$

Integrating once and using $h'(0)=h(0)=f(0)=0$, we get for $t\geq 0$
\begin{equation}
M h'(t)+2h(t)=2f(t)+\int_{s=0}^t\left ( {\psi}(s )
+\phi^2(s)+\phi^3 (s)\right )ds.\label{h**}
\end{equation}
Integrating again, we get
\begin{equation}
h(t)=\frac{1}{M}\int_{s=0}^te^{\frac{2}{M}(s-t)}\left ( 2f(s)+
\int_{r=0}^s({\psi}(r )
+\phi^2(r)+\phi^3 (r))dr\ \right )ds.\label{newha}
\end{equation}
Define
\begin{equation}
\Phi
(t)=\frac{1}{M}\int_{s=0}^te^{\frac{2}{M}(s-t)}\left ( 2f(s)+
\int_{r=0}^s
(\phi^2(r)+\phi^3 (r))dr\right )ds.\label{bigPhia}
\end{equation}
Define the operator $K$ by
\begin{equation}
(Kp)(t)=\frac{1}{M}\int_{s=0}^te^{\frac{2}{M}(s-t)}\int_{r=0}^{s}\int_{w=0}^r\left (
\frac{\partial k}{\partial x}(a^-;0,w)+\frac{\partial k}{\partial x}(a^+;0,w)\right )p(r-w)dwdrds.\label{Ka}
\end{equation}
Thus by \rq{newha},
\begin{equation}
(I-K)h=\Phi ,\label{IKa}
\end{equation}
and formally $h$ is solved by
\begin{equation}
h=\sum_{n=0}^{\infty}K^n\Phi.\label{ha}
\end{equation}
We now prove the convergence of this series. It is easy to show for
$t<T$
$$|K\Phi (t)|\leq t\|\left (
\frac{\partial k}{\partial x}(a^-;0,\cdot)+\frac{\partial k}{\partial x}(a^+;0,\cdot)\right )\|_{L^1(0,t)}\| \Phi\|_{\infty}=Ct\|\Phi \|_{\infty},$$
where $C=\|
\frac{\partial k}{\partial x}(a^-;0,\cdot)+\frac{\partial k}{\partial x}(a^+;0,\cdot)\|_{L^1(0,t)}$.
Here $||*||_{\infty}$ denotes $||*||_{L^{\infty}(0,t)}$.
It follows inductively that
$$|K^n\Phi (t)|\leq \frac{C^n\|\Phi \|_{\infty}t^n}{n!}.$$
This shows that the series converges uniformly on compact sets in $t$, and the solution $h$
is continuous in $t$. That $h\in C^3$ follows from  the following bootstrapping
argument. By  \rq{IKa}, we have
$$h=\Phi +Kh.$$
Since $\Phi \in C^3$ and $h$ is continuous, the smoothing properties of $K$ imply
  $h\in C^1$. Iterating this argument, we get
$h\in C^2$, and then $h\in C^3$.

 That $h(t)=0$ for $t\leq 0$ follows from \rq{newha}. To prove uniqueness of $h$, note that any solution to \rq{Mh} must solve
 \rq{newha}, and hence \rq{IKa}  and \rq{ha} hold as  equations of continuous functions. Since $f$ uniquely determines $\Phi$, \rq{ha} shows that $f$ uniquely
 determines $h$.

We now prove  part B.
By part A,
 for $f\in C^2_*$,
the operator $Sf=h$ is well defined. It is easy to see that $S$ is linear.
Fix $f\in L_*^2$.
We now
estimate $h=Sf$. In what follows, let $C$ be various constants that are independent of $t,f$.
We have
\begin{eqnarray*}
|h(t)| & \leq & \sum_n\frac{C^n||\Phi ||_{\infty}t^n}{n!}\\
& \leq & \sum_n\frac{C^n||f||_{L^2(0,T)}t^n}{n!},
\end{eqnarray*}
so $||h||_{L^2(0,T)}\leq C ||f||_{L^2(0,T)}$.
It follows from \rq{h**} that
\begin{equation}
||h'||_{L^2(0,T)}\leq C ||f||_{L^2(0,T)},\label{hH1}
\end{equation}
 so $h\in H^1_*$. By \rq{hH1}, it follows that $S$ extends to a bounded operator from $L^2_*$ to $H^1_*$.

Next, we prove the invertibility of $S$ in part B. By the Open Mapping Theorem, if suffices to prove S is a bijection. Fix $h\in H^1_*(-\infty ,T)$. In what
follows, terms that are uniquely determined by $h$ will be denoted $F(t)$.
We then rewrite \rq{h**}:
\begin{equation}
F(t)= 2f(t)+
\int_{s=0}^t
\phi^2(s)+\phi^3 (s)ds.\label{solvef}
\end{equation}
We can show that the integral terms on the right hand side can each be expressed in the form $\int_0^tf(s)K(s,t)ds$ with $K$ continuous.
We  prove this for one such integral term, namely
(recalling the $\phi^3$ is given by \rq{phi3a})
$$\int_{s=0}^t
\int_{r=0}^{s}\frac{\partial k} {\partial x}(a^-;0,r)
\int_{w=a}^{s+a-r}k(0^+,a,w)f(s+a-r-w)dwdrds$$
$$=\int_{s=0}^t
\int_{r=0}^{s}\frac{\partial k} {\partial x}(a^-;0,r)
\int_{w=0}^{s-r}k(0^+,a,w+a)f(s-r-w)dwdrds.$$
With the change of variables $\eta =r+w,\ r=r$, the right hand side becomes
\begin{equation}
\int_{s=0}^t
\int_{\eta=0}^{s}f(s-\eta)K_1(\eta )d\eta ds,\label{firstint}
\end{equation}
with
$$K_1(\eta )=\int_{r=0}^{\eta}\frac{\partial k} {\partial x}(a^-;0,r)k(0^+,a,\eta -r+a)dr .$$
With another change variables,  \rq{firstint} equals
$$\int_{r=0}^tf(r)\left ( \int_{s=r}^t K_1(s-r)ds \right ) dr,$$
as desired.
The other integral terms on the right hand side of \rq{solvef} can be treated similarly.
Thus we  can rewrite \rq{solvef} as a Volterra equation of the second kind,
$$F(t)=f(t)+\int_0^tf(s)K(s,t)ds.$$
To prove injectivity of $S$, note that if $h=0$ then $F(t)=0$, and by properties of Volterra equations we conclude $f=0$. Since
$h$ uniquely determines $F(t)$, surjectivity also follows.

The proof of part C is similar to part B, and is left to the reader.

\

{\bf Remarks:}

1- The solution $h$ to  \rq{h1} exists for all $t$, but gives the position of the mass only  for  $t< \min (3a, 2\ell -a)$, ie. until
a reflected wave reaches $x=a$ either from the left or from the right.

\

2- The equations for $u$, \rq{wf1}, \rq{smallx}, \rq{bigx}, all show that for $t<\ell$, the wave propagates at unit speed throughout the interval despite the presence of the mass at $x=a$. The unit speed of propagation will continue to hold for larger times too,
and this will play a key
role in our study of the control problem.

\

3-  It is  important to note that all the calculations above are local, so these results will  extend to the N mass case.

\

4- Lemma \ref{existenceh} together with \rq{smallx} show that a wave reflected off the mass will have the same regularity as the incoming wave.  More precisely, if a wave reaches $x=a$ at time $T$, and if near the mass the function $x\mapsto u(x,T-\epsilon )$ is in $ H^j_{loc}$ for small $\epsilon$, then
$x\mapsto u(x,T+\epsilon )$ is in $ H^j_{loc}$. We now show the same property for a reflection off an endpoint. For simplicity of exposition, we choose
 reflection at $x=\ell $ at time $t=\ell$.
 For $t\leq \ell$ and $x>a$, by \rq{bigx} the function
 \begin{equation}
 w_*^1(x,t):=h(t-x)+\int_{s=x-a}^{t-a}k(a^+;x-a,s)h(t-a-s)ds\label{w1*}
 \end{equation}
 solves the relevant equations of System \rq{pdn1}-\rq{SA}. However, for $t> \ell$,
this expression no longer satisfies $u(\ell, t)=0$. We correct for this by adding the term that models the wave's reflection, denoted $w_*^2(x,t)$, that uniquely solves the
system
\begin{eqnarray*}
v_{tt}-v_{xx}+q(x)v& = & 0, \ t>\ell,\ x<\ell , \\
v(\ell ,t) & = & -w^1_*(\ell, t),\\
v(x,t)& =& 0,\ t<\ell .
\end{eqnarray*}
By the analogue of \rq{reflectg}, we have
\begin{equation}
w_*^2(x,t)=-w_*^1(\ell, t+x-\ell)-\int_{s=\ell-x}^{t-\ell}k(\ell^-;\ell-x,s)w_*^1(\ell, t-s)ds.\label{refw}
\end{equation}
We then have $u=w_*^1+w_*^2$ for $(x,t)$ sufficiently close to $(\ell, \ell)$.
For fixed $t$, it is also clear that $w_*^1\in H^j_*$ implies $w_*^2\in H^j_*$.

\

5- We can solve for $u$ for larger times, as we did in the previous remark, by adding more terms to $u$ and $h$ to account for all the reflections that take place.

\vspace{.5in}

Summarizing, when a wave moves to either an end point or a mass, its reflection will have the same regularity as the incoming wave. Also,
by Lemma \ref{existenceh}, a  wave transmitted across a mass will be one Sobolev order more regular than the incoming wave.
By applying these two principles, the proof of Proposition \ref{Nreg0} for $N=1$ is completed for larger times; the details are left to the reader.
The arguments of this section can also be used to prove Theorem \ref{eur}.

\

\end{subsection}
\begin{subsection}{Regularity and compatibility conditions for $N$ masses}\label{solutionrep}
In this section, we discuss existence, uniqueness, and regularity of solutions to System \rq{A_N11}-\rq{control9}.
We then discuss certain equations that must be  satisfied by functions in the reachable sets; these equations are determined by
the compatibility conditions associated to the wave equation at the masses and at $x=\ell$.
  This will enable us to complete the proof of Theorem \ref{reachable}, and is also necessary for the proof of Theorem \ref{shapecont}.
Consider the System \rq{A_N11}-\rq{control9}, which we rewrite:
\begin{eqnarray}
u_{tt}-u_{xx}+q(x)u& =&0, \ t\in (0,T),\ x\in (0,\ell)\setminus \{ a_j\}_{j=1}^N ,\nonumber\\
u(x,t)& =& 0,\ t\leq 0,\nonumber\\
u(a_j^-,t)& =& u(a_j^+,t)=h_j(t), \ j=1,\ldots,N, \nonumber\\
M_j h_j''(t)& = & u_x(a_j^+,t)-u_x(a_j^-,t),\nonumber\\
u(0,t)& =& f(t),\nonumber\\
  u(\ell ,t)& =& 0.\label{A_N2}
\end{eqnarray}
 We have the following
 extension of Proposition 1.
\begin{prop}\label{Nreg}

For any $T>0$ and $f\in L^2(0,T)$, there exists a unique solution $(u^f,h_1,...,h_N)$ to System \rq{A_N2}. Furthermore,

A) suppose  $T \in ( a_{j},a_{j+1} ] , j = 0, ... N. $ Then
$$ \{ u^f(\cdot,T) : f \in L^2(0,T) \} \subset \{ \phi \in \tilde{W}_0 :
\phi(t) = ... \phi^{(j-1)}(t) = 0,\ \forall t\geq T \}. $$

B) For $ T > \ell , $
 $$ \{ u^f(\cdot,T) : f \in L^2(0,T) \} \subset \{ \phi \in  \tilde{W}_0 :
\phi(l) = 0 \}. $$

C) For any $T>0$,
$$u^f(x,t)\in C([0,T],  \tilde{W}_0) .$$

D) For any $T>0$
$$u^f_t(x,t)\in C  ( [0,T],  \tilde{W}_{-1}) .$$

E) $h_j\in H^j_*$.
\end{prop}
Here $\tilde{W}_i$ were defined in Section 1, and
 for $n\geq 0$ we define,
$$||\phi ||^2_{H^{n}(a_j,a_{j+1})}=||\frac{d^n\phi}{dx^n}||^2_{L^2(a_j,a_{j+1})}+||\phi ||^2_{L^2(a_j,a_{j+1})},$$
while
$$||\phi ||_{\theta'}=\sup_{v\in \theta (0,a_1): ||v||_{H^1}=1}<u,v>,$$
where $<*,*>$ is the duality pairing of $\theta'$ with $\theta$.
We sketch the proof. Uniqueness of $u^f$ follows from Proposition \ref{wavesol}. Existence, along with parts A and B, are essentially proven in
 the previous subsection, where we  demonstrated that a reflected (from a mass or a boundary) wave has the same regularity as
an incident wave, and a transmitted wave is one unit more regular in the Sobolev scale (see, e.g.
representations  \rq{bigx},  \rq{fgh},  \rq{smallx},   and Lemma \ref{existenceh}), and where we also demonstrated unit speed of propagation.
Part C follows easily from the representation of the solution, and part D is  proven using the argument to prove Theorem \ref{eur}  together with
the arguments above.
The details are left to the reader. Part E follows immediately from Lemma \ref{existenceh}.

\

To further clarify the properties of the solutions to our system, we now consider certain equations that must be satisfied
at the masses and at $x=\ell$.
 In what follows, it will be convenient to define
$$L:= -\frac{d^2}{dx^2}+q(x),$$
the differential operator acting on distributions living on $(0,a_1)\cup ...\cup  (a_N, \ell )$.
By Proposition \ref{Nreg},
\begin{equation}
u^f(a_j^-,t)=u^f(a_j^+,t),\label{comp1'}
\end{equation}
 will hold for all $f\in L^2$ and  $j=2,....,N$. Also, assuming $u^f$ is sufficiently regular,
 we have  first order  conditions
\begin{eqnarray}
u^f_x(a_j^-,t) & =  & u^f_x(a_j^+,t)-M_ju^f_{tt}(a_j^+,t)\nonumber \\
& =  &  u^f_x(a_j^+,t)-M_j[u^f_{xx}(a_j^+,t)-q(a_j^+)u^f(a^+,t)]\nonumber\\
& = & u^f_x(a_j^+,t)-M_j(Lu^f)(a_j^+,t).\label{comp2'}
\end{eqnarray}
In particular this will hold for all $f\in L^2$ for $j=3, ..., N$.

In what follows, we write $u$ for $u^f$.
If $u$ is sufficiently regular, then the wave equation also imposes higher order conditions at $x=a_j$. These can be
formulated inductively  from \rq{comp1'} and \rq{comp2'} as follows.
We claim for any positive integer $n$ and $u$ sufficiently regular,
\begin{equation}
L^nu(a_j^-,t)=L^nu(a_j^+,t),\ \forall t>0,\label{evencomp}
\end{equation}
and
\begin{equation}
 \frac{\partial}{\partial x}L^nu(a_j^-,t)=\frac{\partial}{\partial x}L^nu(a_j^+,t)-M_jL^{n+1}u(a_j^+,t)\ \forall t>0.\label{oddcomp}
\end{equation}
We prove \rq{oddcomp} by induction; the proof for \rq{evencomp} is similar but easier.
By \rq{comp2'}, \rq{oddcomp} holds for $n=0$. Assume it holds for some integer $n$. By continuity and the wave equation,
\begin{eqnarray*}
\frac{\partial}{\partial x}L^{n+1}u(a_j^-,t) & = & \lim  _{\ep \to 0^+}\frac{\partial}{\partial x}L^{n}\frac{\partial^2}{\partial t^2}u(a_j-\ep,t)\nonumber\\
& = & \lim  _{\ep \to 0^+}\frac{\partial^2}{\partial t^2}\frac{\partial}{\partial x}L^{n}u(a_j-\ep,t)\nonumber \\
& = & \frac{\partial^2}{\partial t^2}\frac{\partial}{\partial x}L^{n}u(a_j^-,t)\nonumber \\
& = & \frac{\partial^2}{\partial t^2}\big ( \frac{\partial}{\partial x}L^nu(a_j^+,t)-M_jL^{n+1}u(a_j^+,t)  \big ) \\
& = &\lim  _{\ep \to 0^+} \frac{\partial^2}{\partial t^2}\big ( \frac{\partial}{\partial x}L^nu(a_j+\ep,t)-M_jL^{n+1}u(a_j+\ep,t) \big ) \\
& = &\lim  _{\ep \to 0^+} \big ( \frac{\partial}{\partial x}L^n\frac{\partial^2}{\partial t^2}u(a_j+\ep,t)-M_jL^{n+1}
\frac{\partial^2}{\partial t^2}u(a_j+\ep,t)  \big ) \\
& = &\lim  _{\ep \to 0^+} \big ( \frac{\partial}{\partial x}L^{n+1}u(a_j+\ep,t)-M_jL^{n+2}u(a_j+\ep,t) \big ) \\
& = &  \frac{\partial}{\partial x}L^{n+1}u(a_j^+,t)-M_jL^{n+2}u(a_j^+,t) . \ \ \Box
 \end{eqnarray*}

The argument leading to \rq{oddcomp} also leads to
higher order
 boundary conditions  at $x=\ell$, namely
\begin{equation}
L^nu(\ell ,t) =0, \ \forall t>0, \ n=0,1,2,... .\label{compbc}
\end{equation}
Thus for instance,
\begin{equation}
u^f_{xx}(\ell,t) =0 \mbox { and }u^f_{xxxx}(\ell,t)  =   2q'(\ell)u^f_{x}(\ell,t).\label{d4}
\end{equation}

Motivated by \rq{evencomp},\rq{oddcomp}, and \rq{compbc}, we have the following definition.

\noindent{\bf Definition}

Let $j=1,..., N$.
A function $\phi(x)$ satisfies the Condition ${\mathcal C}^{i}$   at $x=a_j$ if
$$ \frac{\partial}{\partial x}L^n\phi (a_j^-)=\frac{\partial}{\partial x}L^n\phi(a_j^+)-M_jL^{n+1}\phi(a_j^+)$$
 is satisfied for $0\leq n\leq \lceil i/2\rceil -2$,
and
$$L^n\phi (a_j^-)=L^n\phi(a_j^+)$$
 is satisfied for $0\leq n\leq \lceil i/2\rceil -1$.
A function $\phi(x)$ satisfies the Condition ${\mathcal C}^{i}$   at $x=\ell$ if
$L^n\phi (\ell) =0$ is satisfied for $0\leq n\leq \lceil i/2\rceil -1$.

For convenience, for $i\leq 0$ we denote the condition
${\mathcal C}^{i}$ at $x=a_j$ to be a vacuous condition.

A function satisfies Condition $\mC^{i}_*$ if it satisfies Condition ${\mathcal C}^{j-1+i}$ at $x=a_j$ for all $j=1,..., N+1$.

\

Also, $u_t$ must satisfy compatibility conditions.
For instance, for $f\in L^2_*$, it follows from Proposition \ref{Nreg} that (compare with \rq{comp1'},\rq{comp2'})
 $$u^f_t(a_j^-,t)=u_t^f(a_j^+,t),\mbox{ if } j\geq 3\mbox{ and }$$
 $$u^f_{xt}(a_j^-,t)= u^f_{xt}(a_j^+,t)-M_j[u^f_{xxt}(a_j^+,t)-q(a_j^+)u^f_t(a^+,t)], \mbox{ if } j\geq 4.$$
It follows   that  for fixed $T>0$,  if $u^f(x,T)$ satisfies Condition ${\mathcal C}_*^{j}$,
then  $u_t^f(x,T)$ satisfies Condition ${\mathcal C}_*^{j-1}$.

We now define for $i=0,-1$ the spaces
$$ W_i:=\{ \phi \in \tilde{W_i} :
\ \phi \mbox{ satisfies Condition $\mC^{i}_*$  }\} . $$	
Then $W_i$ are  real Hilbert spaces with norms
$$
||u||^2_{W_0}=\sum_{j=0}^N||u||^2_{H^{j}(a_j,a_{j+1})}+\sum_{j=1}^NM_ju(a_j^+)^2, $$
$$||u||_{W_{-1}}^2=||u||_{\theta '}^2+\sum_{j=1}^N||u||^2_{H^{j-1}(a_j,a_{j+1})}+\sum_{j=2}^NM_ju(a_j^+)^2 .\label{hilbert1}
$$

Theorem 1 now follows from Proposition \ref{Nreg} and our definitions of $W_0,W_{-1}$.

\end{subsection}

\end{section}
\begin{section}{Shape and velocity controllability results}
The goal of this section is to prove Theorem \ref{shapecont}.
To prove  this theorem for $T>\ell$,
we require a preliminary result that we formulate now.
Consider
the vibrating string with no mass:
\begin{eqnarray}
u_{tt}-u_{xx}+q(x)u& =&0, \ t\in (0,T),\ x\in (0,\ell),\label{mtl}\\
u(x,t)& =& 0,\ t\leq 0,\\
u(0,t)& =& g(t),\\
  u(\ell ,t)& =& 0.\label{bcnhN}
\end{eqnarray}
Here we will assume $q\in C^{\max(0,N-2)}[0,\ell ]$.
Let
$${\cal A}_N=\{ (\phi (\ell ), \phi^{(1)}(\ell) ,...,\phi^{(N-1)}(\ell))\in \bR^N: \ \phi \in H^N(0,\ell)  \mbox{ and } \phi \mbox{ satisfies Condition }{\cal C}^N\mbox{ at }x=\ell \} $$
$$=\{ (d_0,...,d_{N-1})\in \bR^N: \exists \mbox{ $T, g$ so that $u=u^g$ solves }\frac{\partial^ju}{\partial x^j}(\ell ,T)=d_j, \ j=0,..., N-1\} .$$
\begin{lemma}\label{pointwise}
	
Let $T>\ell$, and let
 $\delta =T-\ell$.

\noindent i) Let ${\bf d}=(d_0,..., d_{N-1})\in {\cal A}_N$.
Then there exist
 $g\in H_0^{N} (0, 2\delta )$ and constant $C$ independent of ${\bf d}$
 such that the solution $u=u^g$ to the system \rq{mtl}-\rq{bcnhN} satisfies
\begin{equation}
 \frac{\partial^ju}{\partial x^j}(\ell ,T)=d_j, \ j=0,..., N-1, \label{onto}
 \end{equation}
 and
 \begin{equation}
 ||g||_{H^N}^2\leq C\sum_{j=0}^{N-1}|d_j|^2.\label{gbound}
 \end{equation}

\noindent ii) Let ${\bf d}=(d_0,..., d_{N-2})\in {\cal A}_{N-1}$.
 Then there exist $g\in H_0^{N} (0, 2\delta )$ and constant $C$ independent of $\phi$ such that the solution $u$ to the system \rq{bcnhN} satisfies
 $$\frac{\partial^{j+1}u}{\partial x^j\partial t}(\ell ,T)=d_j, \ j=0,..., N-2$$
 and
 $$||g||_{H^{N}}^2\leq C\sum_{j=0}^{N-2}|d_j|^2.$$
\end{lemma}
Proof:
We prove part i; the simple adaptation for part ii is left to the reader.
We define the following linear functionals on $H_0^N(0,2\delta)$:
$$\tau_j (g):=\frac{\partial ^ju}{\partial x^j}(\ell, T).$$
 Consider now the bounded operator $A_1: H_0^{N}(0,2\delta )\ra \bR^N$ given by
$$A_1(g)=(\tau_0(g),..., \tau_{N-1}(g)).$$

Proving \rq{onto} is equivalent to proving $A_1$ maps onto ${\cal A}_N$.

Claim 1: $dim \ {\cal A}_N \leq  r:=\lfloor N/2\rfloor$.

Proof of claim: recall that  Condition ${\mathcal C}^{N}$ at $x=\ell$ imposes a set of apriori equations on
$\{ \phi (\ell ), ...,\phi^{(N-1)}(\ell)\}$, determined by \rq{mtl} in tandem with the condition $u(\ell, t)=0$.
We will now present $r$ such equations.
Let $w$ be any solution
to \rq{mtl} with $w(\ell, t)=0$. Assuming $w$ is sufficiently regular,
we have for $x<\ell $ that $L^nw(x,t)=\frac{\partial^{2n}w}{\partial t^{2n}}(x,t)$, from which we conclude
$$\frac{\partial^{2n}w}{\partial x^{2n}}(\ell ,t)=\sum_{j=0}^{2n-2}c_{j,n}\frac{\partial^{j}w}{\partial x^{j}}(\ell ,t),$$
where the constants
$c_{j,n}$ depend of $q$.
Thus
$$\phi^{(2n)}(\ell) =\sum_{j=0}^{2n-2}c_{j,n}\phi^{(j)}(\ell )$$
for any $\phi\in H^N(0,\ell )$ satisfying Condition ${\cal C}^N$ at $x=\ell $.
Such an equation holds for each $n$ such that $0\leq 2n\leq N-1$, and there are $\lceil N/2\rceil $ such $n$, so the claim follows.

%involving boundary conditions at $x=ell$. Each such equation imposes a linear independence relation among $\{ \tau_0,...,\tau_{N-1}\}$.
%or instance,  for any $g$, since $u(\ell,t)=u_{xx}(\ell,t)=0$, we have $\tau_0=\tau_2=0$, and by \rq{d4} we have  $\tau_4=-2q'(\ell)\tau_1$.
% We will show that for each $n$ such that $2n\leq N_1$,
% \begin{equation}
% \tau_{2n}=\sum_{j=0}^{2n-2}c_j^n\tau_j,  \mbox{ for constants }c_j^n.\label{tauj}
% \end{equation}
%The claim follows immediately.

Claim 2: the rank of $A_1$ is  $r$.

Proof of claim:  for any smooth function $g(t)$,
for $(x,t)$ sufficiently close to $(\ell, \ell)$, we have by analogues of \rq{w1*} and \rq{refw} that
$u(x,t)=w_*^1(x,t)+w_*^2(x,t)$ with
$$w^1_*(x,t)=g(t-x)+\int_{s=x}^tk(0^+;x,s)g(t-s)ds$$
and
$$w^2_*(x,t)=-g(t+x-2\ell)-\int_{s=\ell}^{t+x-\ell}k(0^+;\ell,s)g(t+x-\ell-s)ds$$
$$-\int_{s=\ell-x}^{t-\ell}k(\ell;\ell-x,s)[g(t-s-\ell)+\int_{r=\ell}^{t-s}k(0^+;\ell,r)g(t-s-r)dr]ds.$$
A tedious but straightforward calculation shows that
$$\frac{\partial ^ju}{\partial x^j}(\ell ,T)=\sum_{k=0}^jC_{j,k}g^{(k)}(T-\ell)+\int_{s=0}^{T-\ell}K_j(s)g(T-\ell -s)ds,\ j=0,..., N-1,$$
with $C_{j,k}$ various constants and $K_j(s)\in C[0,T-\ell ]$, and
furthermore
$$C_{j,j}=-2 \mbox{ if $j$ is odd }, \ C_{j,j}=0 \mbox{ if $j$ is even}.$$
Let $j_0$ be the largest  integer such that $2j_0-1\leq N-1$.
We show  the set ${\cal T}:= \{ \tau_{1},\tau_3,...,\tau_{2j_0-1}\}$ is linearly independent.
Let $\gamma >0$. We can choose $g\in H^N(0,2\delta)$ satisfying:
$$g^{(2j_0-1)}(\delta) =1, \ g^{(j)}(\delta) =0, j=0,...,2j_0-2,\mbox{ and }\int_{s=0}^{T-\ell}K_j(s)g(T-\ell -s)ds<\ga,\ j=0,...,N-1.$$
Let $\ep>0$.
Then choosing $\ga$ sufficiently small,
$$|\tau_{2j_0-1}(g)|>1-\ep \mbox{ and }|\tau_{j} (g)|<\ep \mbox{ for } j=0,..., 2j_0-2.$$
This proves $\tau_{2j_0-1}$ is linearly independent of the span of $\tau_0,...,\tau_{2j_0-2}$.
Iterating this argument, it follows that ${\cal T}$ is linearly independent,
so the rank of $A_1$ is at least $r$. Claim 2 now follows from Claim 1.

%Since $u(\ell,t)=u_{xx}(\ell,t)=0$, we have $\tau_0=\tau_2=0$, and by \rq{d4} we have  $\tau_4=-2q'(\ell)\tau_1$.
%More generally,
% linear dependence relations for a subset of $\{ \tau_j\} $, which are determined by the equation $u_{tt}-u_{xx}+qu=0$ in tandem with $u(\ell ,t)=0$, are in one to one correspondence with  the boundary conditions associated Condition ${\mathcal C}^{N}$ at
% $x=\ell$.
% Let $D$ be the dimension of the $V:=\mbox{span}(\tau_0,...,\tau_{N-1}).$ Let $\sigma: \{ 1,..., D\} \ra \{ 0,..., N-1\}$ be an injection such
% that $\{ \tau_{\sigma(j)}\}_{j=1}^D$ is a basis of $V$. We claim the mapping $\tau: C_0^{\infty}(0,2\delta ) \ra \bR^D$ given by
%$\tau (g)=(\tau_{\sigma (1)}(g), ...,\tau_{\sigma(D)}(g))^T$ is onto.  In fact, if not, then there would exist a non-zero vector ${\bf %c}=(c_1,...,c_D)^T\in \bR^D$
% with $$0={\bf c}\cdot \tau (g)=\sum_{j=1}^Dc_j\tau_{\sigma (j)}(g),\ \forall g\in C_0^{\infty}(0,2\delta ) .$$
% This would contradict the linear independence of $\{ \tau_{\sigma(j)}\}_{j=1}^D$.
% Finally, given $\phi$ satisfying the hypotheses of the lemma, we find $g\in C_0^{\infty}(0,2\delta )$ solving
%$$\tau_{\sigma (j)}(g)=\phi^{\sigma (j)}(\ell),\ j=1,..., D.$$
% The remaining boundary conditions will automatically be satisfied because of the linear dependence relations.

It follows from Claim 1 and Claim 2 that $A_1$ maps $H^N_0(0,2\delta )$ onto ${\cal A}_N$.

We now prove \rq{gbound}.
Let $g_1,...,g_r$ be linearly independent functions in $H_0^{N}(0,2\delta)$  such that $\{  A_1g_1,...,A_1g_r\}$ span ${\cal A}_N$. Let ${\bf b}=(b_1,...,b_r)$, and define
the operator $A_2$ by $A_2({\bf b})=b_1g_1+...+b_rg_r$, and let  $A=A_1A_2$.
%Define $B:\bR^r\ra \bR^N$ by $B(b_1,...,b_r)=A(b_1g_1+...+b_rg_r).$
Then clearly $A:\bR^r\ra \bR^N$ is an injection, and setting $g=b_1g_1+...+b_rg_r$ and $u=u^g$, there exist positive constants $C_1,C_2$ such that
$$\sum_{k=0}^{N-1}  |\frac{\partial^{k}u}{\partial x^k}(\ell, T)|^2 =||A_1g||_{{\bf \bR}^N}^2=||A{\bf b}||^2_{{\bf \bR}^N}\geq C_1 \sum_1^r|b_j|^2\geq C_2||g||^2_{H^N(0,2\delta)}.$$
The proof is complete.

\vspace{.5in}

We now prove part A of Theorem \ref{shapecont}.
The proof of part B, which is similar, is left to the reader.
Recall
$W_i^T=\{ \phi \in W_i: \phi (x)=0\mbox{ for }x\geq T\}$.
Thus $\phi \in W_0^{\ell}$ implies
$\phi^{(j)}(\ell)=0$ for $j=0,..., N-1,$ whereas
 $\phi \in W_0^{}$  implies only $\phi (\ell)=0$.

We need to prove the following statements:

\noindent{\bf I}.  Suppose  $T \leq \ell .$  For any  $\phi \in W^T_0$, there exists a unique $f\in L^2(0,T)$
such that $u^f(x,T)=\phi(x)$. Furthermore,
\begin{equation}
\| u^f(\cdot,T)\|_{W_0} \asymp \| f\|_{L^2(0,T)}.
\label{lb'}
\end{equation}
\noindent{\bf II}. Let $ T > \ell . $
For any  $\phi \in W_0$, there exists a  $f\in L^2(0,T)$
such that $u^f(x,T)=\phi(x)$ and \rq{lb'}  holds.

\

Before proving theorem in the general case,
it will be instructive to prove part I in the simpler case with no mass. Thus, let $T\leq \ell$,
and let
$\phi \in L^2(0,\ell )$ satisfy $\phi (x)=0$ for $x>T$. We wish to find $f\in L^2(0,T)$ such that $u=u^f$ solving
\begin{eqnarray*}
u_{tt}-u_{xx}+q(x)u& =&0, \ t\in (0,T),\ x\in (0,\ell),\nonumber\\
u(x,t)& =& 0,\ t\leq 0,\nonumber\\
u(0,t)& =& f(t),\nonumber\\
  u(\ell ,t)& =& 0
\end{eqnarray*}
satisfies
$$  u(x,T)= \phi (x), \ \forall x.$$
By Proposition \ref{wavesol},
\begin{equation}
\phi(x)=
u(x,T)=
f(T-x)+\int_{s=x}^Tk(0^+;x,s)f(T-s)ds.\label{babycase}
\end{equation}
This is a Volterra equation of the second kind, so there exists a unique solution $f\in L^2(0,T)$, thus proving shape controllability  in this simple case.

In the case of attached masses, the function $u(x,T)$ does not satisfy an equation like \rq{babycase} for all $x$ because there will be addition terms
on the right hand side of \rq{babycase}  due to reflections off the masses. Nevertheless, the argument above can be adapted to the system with masses because there
will also be an interval near the wavefront, at $x=T$, where there will be no reflected waves, so a representation like \rq{babycase} holds.

\

\noindent{\bf  Proof of part I.}
Fix $T\leq \ell$ and suppose  $j$ satisfies $T\in (a_j,a_{j+1}]$.
We will prove that for any $\phi(x) \in W^T_0,$ the equation
\beq \lab{cty}
 u^f(x,T) = \phi(x)
\eeq
has a unique solution  $ f \in L^2(0,T)  $.
  We set $\La = 2 \min|a_{i+1} - a_i|, \;  i = 0, \ldots , N,$ and will solve the equation \rq{cty}
by steps of the length at most $\La.$ This means that we will move by such steps along the $x-$axis from the right to the left starting at the point $x=T$.

{\bf Step 1.} We solve the equation \rq{cty} for the set of $x \in (a_j,T)$ within distance $\Lambda$ of the the wavefront $x=T$.
On this interval at time $t=T$,
 there can be no reflected terms because of unit speed of propagation.
We consider separately the  two possible cases: (a)  $T-\La \geq a_j, \ $ (b)  $T-\La <a_j.$

Case a:  Because there are no reflected terms on the interval
 $x>T-\La$, arguing as in \rq{bigx} we have  that the equation \rq{cty} reduces to
\beq \lab{ctj}
\phi(x) = h_j(T-x) + \int_{x-a_j}^{T-a_j} k(a_j^+;x-a_j,s) h_j(T-a_j-s) ds,
\eeq
with $x\in (T-\Lambda, T]$.
This is a Volterra equation of the second kind, hence $\phi$ uniquely determines the function $h_j(t)$ on the interval $[0, \La ),$ and its regularity is the same as the regularity of
$\phi.$ Using Lemma 1 we can conclude there exists a unique $f\in L^2(0,\La )$ such that
 $f(t)=(S^{-j}h_j)(t)$ for $t \in [0,\La).$

Case b: We use the same argument as in part a over the shorter interval $(a_j, T]$. Thus we
  first use \rq{ctj} on the interval $x \in (a_j, T]$ to find  $h_j(t)$ on the interval $[0, T-a_j).$ Then we find
$f(t)= (S^{-j}h_j)(t)$ on the same time interval.
%For the reminder of this proof, it will be convenient to denote $\ell_j=a_j-a_{j-1}$.
Define $c_1$ by $c_1=\La$ in case a, and $c_1=T-a_j$ in case b.

{\bf Step 2}
 Define $f_1$ by
$$f_1(t)=\left \{
\begin{array}{cc}
f(t), & t< c_1 ,\\
0, & t\in [c_1,T] .
\end{array}
\right .
$$
Then by construction, $\phi (x)=u^{f_1}(x,T)$ for $x>T-c_1$, and hence
$\phi_2(x):=\phi (x)-u^{f_1}(x,T) \in W_0^{T}$ is supported on $[0,T-c_1]$, and hence is in $W_0^{T-c_1}$.
% Assume for the moment that $T-c_1>0$.
If $T-c_1>a_j$,  then  we can repeat the argument in Step 1. If $T-c_1\leq a_j$, then the same argument as in Step 1 can be
carried out in the interval $(a_{j-1},a_j]$. In both cases,  we find $c_2$ along with the unique $f_2\in L^2(0,T)$ supported in $[c_1 ,c_1+c_2 ]$
such that $u^{f_2}(x,T)=\phi_2 (x)$ for $x\in [T-c_1-c_2, T-c_1 ]$.

Thus $u^{f_1+f_2}(x,T)=u^{f_1}(x,T)+u^{f_2}(x,T)=\phi (x)$ for $x\in [T-c_1-c_2, T ]$

{\bf Step 3}
We repeat the arguments of step 1 and  step 2 as often as necessary, thus solving for $f$.

\

We now prove \rq{lb'}.  Suppose
there exists a sequence $\{ f^n\} $, with $f\in {L^2(0,T)}$ and $ \{ \phi^n=u^{f^n}(x,T)\}$, such that $\| \phi^n\|_{W_0}\to 0$. We will prove
$\| f^n\|_{L^2(0,T)}\to 0$.

 Consider Step 1. Suppose $h_j^n(t)$ solves  \rq{ctj}, i.e.
\beq \lab{ctj'}
\phi^n(x) = h_j^n(T-x) + \int_{x-a_j}^{T-a_j} k(a_j^+;x-a_j,s) h_j^n(T-a_j-s) ds, \ T-c_1 < x < T .
\eeq
By assumption, we have
$$\| \frac{d^m}{dx^m}\phi^n\|_{L^2(T-c_1,T)} \to 0,\ \ m=0,...,j.$$
We apply this with  $m=0$ to \rq{ctj'}, which is a Volterra equation of the second kind.
Arguing as in the proof of Lemma \ref{existenceh} part B,
we deduce $\| h_j^n\|_{L^2(0,c_1)}\to 0$.
An inductive argument applied to $x$-derivatives of \rq{ctj'} then easily gives
$$\| \frac{d^m}{dt^m}h_j^n\|_{L^2(0,c_1 )} \to 0,\ \ m=0,...,j.$$
It then follows by Lemma \ref{existenceh} that $\| f^n\|_{L^2(0,c_1 )}\to 0.$

We now use the ideas in Step 2 to finish the proof. Specifically,
let
$$f^n_1(t)=\left \{
\begin{array}{cc}
f^n(t), & t< c_1 ,\\
0, & t\geq c_1 .
\end{array}
\right .
$$
Then ${\phi}_2^n(x):=\phi^n (x)-u^{f_1^n}(x,T) \in W^T_0$ is supported in $[0,T-c_1 ]$, and
$\| {\phi}^n_2\|_{W_0}\to 0$. Thus we can repeat the argument above on ${\phi}_2^n$.
Now iterating this argument as many times as necessary, as in Step 2, we get $\| f^n\|_{L^2(0,T)}\to 0$. This proves \rq{lb'}.

\noindent{\bf Proof of part  II.}\

 Let $\phi \in W_0$.
Let $T=\ell +\delta$, where we can assume without loss of generality
$\delta \in (0, \La /4)$. Recall $\ell_N=\ell-a_N$.
% Let $\rho(x)\in C^{\infty}(a_N,\ell)$ such that the support of $\rho$
%is $[\ell-2\delta,\ell ]$, and $\rho =1 $ on $[\ell-\delta, \ell ]$.
Consider the system
\begin{eqnarray*}
\tilde{u}_{tt}-\tilde{u}_{xx}+q(x)\tilde{u}& =&0, \ t\in (0,\ell_N+\delta),\ x\in (a_N,\ell),\\
\tilde{u}(x,t)& =& 0,\ t\leq 0,\\
\tilde{u}(a_N,t)& =& g(t),\\
\tilde{u}(\ell ,t)& =& 0.\\
\end{eqnarray*}
By Lemma \ref{pointwise},
 there exists  $g(t)\in H_0^{N}(0,2\delta )$ such that $\tilde{u}=\tilde{u}^g$ satisfies
$$\frac{\partial^j\tilde{u}}{\partial x^j}(\ell,\ell_N+\delta )=\phi^{(j)} (\ell), \ j=0,..., N-1.$$
Since $g$ is in the range of $S^N$, we can set $f_g=S^{-N}g\in L^2(0,2\delta )$, so that
 $\phi(x) -u^{f_g}(x,T)\in W_0^{\ell}.$
%Because $\tilde{u}^g(x,T-a_N)$ is differentiable to order $N-1$, it
%it is not hard to see that $g^{(n)}(0^+)=0$ for $n=0,..., N-1$. In particular, $g$ is in the range of $S^N$.
%We claim there exists a control $g(t)\in H_*^N(0,T-a_N )$ such that
%$\tilde{u}^g(x,T-a_N)=\rho (x)\phi (x).$ Now choose a cutoff function $\xi\in C^{\infty}([0,T-a_N])$
%such that $\xi=1$ on $[0,2\delta ]$,  $\xi(t)=0$ for $t>3\delta$.
%Then $\tilde{u}^{\xi g}(x,T)=\rho (x)\phi(x)$ for $x\in [\ell -\delta ,\ell]$.
%, and $u^{\xi g}(x,T)=0$ for $x\in [0,\ell -2\delta ]$.
%Let $f_g=S^{-N}(\xi g)$, so $f_g\in L^2$ is supported in $[0, T-a_N]$.
%We have $u^{f_g}(x,T)\in W$, and $\phi(x) -u^{f_g}(x,T)=0$ for $x>\ell -\delta$. This means $\phi(x) -u^{f_g}(x,T)\in W^{T-2\delta}$.
 Thus
 we can apply the part i to find $f\in L^2(0 , T)$, supported in $(\delta,T)$, solving $u^f(x,T)=\phi(x) -u^{f_g}(x,T)$. Then
$f+f_g$ is the desired control.

To prove \rq{lb'} in this case, suppose
there exists a sequence $ \{ \phi^n\}$ such that $\| \phi^n\|_{W_0}\to 0$, and suppose we use the procedure from the previous paragraph to find
$\{ f^n\} $, with $f^n\in {L^2(0,T)}$ such that $\phi^n=u^{f^n}(x,T)$.
 We will prove
$\| f^n\|_{L^2(0,T)}\to 0$.

Since $\| \phi^n\|_{W_0}\to 0$, it follows that the restriction of $\phi^n$ to the interval $(a_N,\ell)$ tends
to zero in $H^N$ topology, and  hence
$$\sum_{j=0}^{N-1} |(\phi^n)^{(j)}(\ell)|^2\to 0.$$
By Lemma \ref{pointwise} part i, the associated sequence of functions $\{ g^n\} $ will converge in zero in $H^N(0,2\delta)$.
Since $g^n$ are in the range of $S^N$, we set $f_{g^n}:=S^{-N}g^n\in L^2(0,2\delta )$.
Then
$$\phi^n(x)-u^{f_{g^n}}(x,T)=u^{f^n-f_{g^n}}(x,T)
\in W_0^{\ell}.$$
Since $f_{g^n}$ converges to zero in $L^2(0,2\delta )$, it follows that
$$||\phi^n(x)-u^{f_{g^n}}(x,T)||_{W_0}\to 0.$$
By part I, we conclude
$||{f^n-f_{g^n}}||_{L^2(0,T)}\to 0$, so $\| f^n\|_{L^2(0,T)}\to 0$. The proof of part II is complete.

\

{\bf Remark:} Using the ideas in this section, we can prove the following shape controllability result. Assume the hypotheses of Lemma \ref{pointwise} part i. Then there exists $f\in H_*^N(0,T)$ satisfying $u^f(x,T)=\phi(x)$. Theorems proving that highly regular target functions can be generated
by highly regular controls have been proven in the  context of $full$ controllability in a number of works, see \cite{EZ} and references therein.

\end{section}

\begin{section}{Spectrum and Riesz bases}
We begin this section, in Subsection 4.1, with a study of the spectrum of the Sturm-Liouville problem associated to System \rq{A_N11}-\rq{control9}.
Let $\{ \la_n\}$ be the associated frequency spectrum. In section 4.2, we
use our results about the asymptotics of the frequencies and
  the eigenfunctions to construct a Riesz basis of $L^2(0,2\ell )$ consisting of E.D.D. of
  the   family  $\{ e^{i\la_nt}\}$. This Riesz basis will be one of the ingredients in the proof of Theorem \ref{mainthm}.
Finally, in Section 4.3, we make some more remarks about the asymptotics of the eigenfunctions.

\begin{subsection}{Spectral theory of system}

Consider the eigenvalue problem which arises by applying separation of variables to System \rq{A_N11}-\rq{control9}:
\begin{eqnarray}
-\phi''(x)+q(x)\phi(x) & = & (\la)^2\phi(x),\ x\in (0,\ell )\setminus \{a_j\}_{j=1}^N,\nonumber \\
\phi(0)& = & \phi(\ell)  =  0,\nonumber\\
\phi(a_j^-)& = & \phi(a_j^+),\nonumber\\
\phi'(a_j^+) & = & \phi'(a^-_j)-M_j(\la)^2\phi(a_j),\ j=1,\ldots , N.\label{Nmassq}
\end{eqnarray}
Let $\{ (\la_n)^2\}_{n=1}^{\infty}$ be the set of eigenvalues of System \rq{Nmassq},
listed in increasing order.
 It is easy to show that the eigenvalues are simple. In fact, suppose $\varphi_1,\varphi_2$ solve
the system \rq{Nmassq} for the same $\la$. Then by scaling, we can also assume $\varphi_1'(0)=\varphi_2'(0)$. Thus $v=\varphi_1-\varphi_2$ satisfies the
$-v''(x)+q(x)v(x) = (\la )^2v(x)$ on $(0,a_1)$ with $v(0)=v'(0)=0$. By standard theory of differential equations, we conclude
$v(x)=0$ on $(0,a_1)$, and hence $\varphi_1(x)=\varphi_2(x)$ on $[0,a_1]$. Applying the same argument to $\varphi_1,\varphi_2$ on the subinterval
$[a_1,\ell]$, we deduce $\varphi_1=\varphi_2$ on $[a_1,a_2]$. Iterating, we get $\varphi_1=\varphi_2$.

The set $\{ (\la_n)^2\}$ can be realized as the spectrum of a self-adjoint operator as follows.
First we introduce
$$\mathcal{H}_0:=L^2_M(0,\ell)=\{ u\in L^2(0,\ell ): \ u(a_j)\in \bR\},$$
 where the norm in $L^2_M(0,\ell)$ is defined as
\begin{equation}
\int_0^{\ell} |u(x)|^2dx + \sum_{j=1}^N M_j|u(a_j)|^2. \label{L2M}
\end{equation}
Denote by $<*,*>_M$ the associated inner product.
This space is canonically isomorphic to $L^2(0,\ell )\oplus \bR^N$.
We define quadratic form
$$Q(u,v)=\sum_{j=0}^N\int_{a_j}^{a_{j+1}}u'(x)v'(x)+q(x)u(x)v(x)dx,$$
with domain
$${\cal Q}=\{ u\in L^2_M(0,\ell ): \ u|_{(a_j,a_{j+1})}\in H^1(a_j,a_{j+1}),\ u(a_j^-)=u(a_j)=u(a_j^+) \ \forall j,\mbox{ and } u(0)=u(\ell)=0\}.$$
Associated with this semibounded, closed quadratic form is the self adjoint operator $A$, with operator
domain $$D(A)=\{ u\in {\cal Q}: \ Au\in L^2_M(0,\ell )\} .$$
Then for $u\in D(A)$,
$$Au(x)=
\left \{
\begin{array}{cc}
-u''(x)+q(x)u(x), & x\neq a_j,\ j=1,...,N,\\
\frac{1}{M_j}(u'(a_j^-)-u'(a_j^+)),& x=a_j,\ j=1,...,N.
\end{array}
\right .
$$
Letting $\{ \varphi_n\}$ be the set of normalized eigenfunctions of $A$, by simplicity of the spectrum and the self-adjointness of $A$ we have that this set is orthonormal with respect to $<*,*>_M$. It is easy to check that $\varphi_n$ will solve the system \rq{Nmassq}.

% We use the spectral representation to create a scale of Sobolev spaces:
%$${\mathcal H}_p=\{ u(x)=\sum_{n=0}^{\infty} a_n \phi_n(x) \;:\; ||u||_p^2 = \sum_{n=0}^{\infty} |a_n|^2  |\la_n|^{2p} < \infty \}, \ p \in \bR.$$
%Thus ${\mathcal H}_2=D(A)$. For $p<0$, we require the assumption that zero is not an eigenvalue of $A$.
%Associated to these spaces are various symmetric compatibility conditions at $x=a_j$. To give an example,  any
%eigenfunction $\phi $ is in ${\mathcal H}_n$ for all $n$,  so that  for each $j$ we have
%$$\frac{1}{M_j}(\phi'(a_j^-)-\phi'(a_j^+))=\la^2 \phi(a_j)= \phi''(a_j^-)=\phi''(a_j^+).$$

\

It is necessary to study in detail the eigenfunctions in the case $q=0$.
Thus
consider the  initial value problem
\begin{eqnarray}
-\phi''(x) & = & (\la )^2\phi(x),\ x\in (0,\ell )\setminus \{a_j\}_{j=1}^N,\nonumber \\
\phi(0)& = &  0, \ \phi'(0) = \la, \nonumber\\
\phi(a_j^-)& =& \phi(a_j^+),\nonumber\\
\phi'(a_j^+) & = & \phi'(a^-_j)-M_j(\la )^2\phi(a_j),\ j=1,\ldots , N.\label{Nmass}
\end{eqnarray}
We solve for $\phi=\phi(x,\la)$ by the following procedure. On the interval
$(0,a_1)$ we set
\begin{equation}
\phi(x,\la)=\sin (\la x),\ x\in (0,a_1).\label{0a_1}
\end{equation}
 We then obtain $\phi(a_1,\la)$ and $ \phi'(a_1^+,\la)$
from equation \rq{Nmass}. Then
$$
\phi(x,\la)=\phi(a_1,\la)\cos (\la (x-a_1))+\frac{\phi'(a_1^+,\la)}{\la}\sin (\la (x-a_1)),\ x\in (a_1,a_2).
$$
Clearly, we can then iteratively solve for $\phi$ on $(a_j,a_{j+1})$ for each $j=0,..., N$.
We will require the following technical lemma. Recall $\ell_j=a_{j+1}-a_j$.
\begin{lemma}\label{q0asy}
%Set $M_0=1$.

\

i) For each $j$, $j=1,\ldots , N+1$,
\begin{equation}
\phi(a_j,\la)=\sum_{n=0}^{j-1}b_n^j(\la)\la^n,\label{bn}
\end{equation}
where for each $j,n$, the $b_n^j(\la )$ is a   linear combination, with coefficients independent of $\la$, of products of the form
\begin{equation}
\prod_{k=0}^{j-1}T_k(\la \ell_k),\label{Tk}
\end{equation}
where  $T_k(z)$ equals either $\sin (z)$ or $\cos (z)$, and at least $(n+1)$ of the $T_k(z)$ equal $\sin (z)$.
Furthermore
$$
b^j_{j-1}(\la)=(-1)^{j-1}\left (\prod_{k=1}^{j-1}M_k\right ) \prod_{k=0}^{j-1}\sin (\la \ell_k).
$$
Here we use the convention $\prod_1^0=1$.

ii)  For each $j$, $j=1,\ldots , N$,
 $$\phi'(a_j^+,\la)/\la =\sum_{n=0}^{j}\tilde{b}^j_n(\la)\la^n,$$
where for each $j,n$, the $\tilde{b}_n(\la )$ is a  linear combination, with coefficients independent of $\la$, of products of the form
$$\prod_{k=0}^j \tilde{T}_k(\la \ell_k),$$
where  $\tilde{T}_k(z)$ equals either $\sin (z)$ or $\cos (z)$,  and at least $n$ of the $\tilde{T}_k(z)$ equal $\sin (z)$.
Furthermore,
$$
\tilde{b}_j^j(\la) =(-1)^{j}\prod_{k=0}^{j-1}M_k \sin (\la \ell_k).
$$

iii)  For $\la$ close to 0, $\phi(a_j,\la)=\sin (\la a_j)+O(\la^3)$ and $\phi'(a_j^+,\la)=\la \cos (\la a_{j})+O(\la^3)$.

iv) $\phi(a_j,-\la)=-\phi(a_j,\la)$ for all $\la$.

\end{lemma}
{\bf Example} We will show $\phi(a_1,\la)=\sin (\la a_1)$, so $b_0^1=1$, and
\begin{equation}
\ \phi(a_2,\la) = -\la M_1\sin (\la \ell_0)\sin (\la \ell_1)
+[\sin (\la \ell_0)\cos (\la \ell_1)+\cos (\la \ell_0)\sin (\la \ell_1)],\label{a2}
\end{equation}
so
$$\  b_1^2= -M_1\sin (\la \ell_0)\sin (\la \ell_1),\  b_0^2= \sin (\la \ell_0)\cos (\la \ell_1)+\cos (\la \ell_0)\sin (\la \ell_1).$$

\noindent{\bf Proof of Lemma:}
For $j=1$, by \rq{0a_1} and \rq{Nmass} we have
\begin{equation}
\phi(a_1,\la)=\sin (\la a_1),\ \ \phi'(a_1^+,\la)/\la =\cos (\la a_1)-M_1\la \sin (\la a_1),\label{j=1}
\end{equation}
so the lemma holds in this case. We now suppose $j=2$. On the interval $(a_1,a_2)$, we have
\begin{equation}
\phi(x,\la)=\phi (a_1,\la)\cos (\la(x-a_1))+\frac{\phi '(a_1^+,\la)}{\la }\sin (\la (x-a_1)).\label{soln12}
\end{equation}
Hence \rq{a2} follows from  \rq{j=1},
and hence clearly part i) is satisfied in this case. Differentiating \rq{soln12} and applying \rq{Nmass}, we get
\begin{eqnarray*}
\frac{\phi '(a_2^+,\la)}{\la}& =& \la^2M_1M_2\sin(\la \ell_0)\sin(\la \ell_1)+ \la \big ( M_2\sin(\la \ell_0)\cos(\la \ell_1)+M_2\cos(\la \ell_0)\sin(\la \ell_1)\\
& &  -M_1\sin(\la \ell_0)\cos(\la \ell_1)\big )-\sin(\la \ell_0)\sin(\la \ell_1)-\cos(\la \ell_0)\cos(\la \ell_1),
\end{eqnarray*}
so part ii) holds for $j=2$.
It is easy to complete the proof of i), ii) by induction on $j$ using
\begin{equation}
\phi(x,\la)=\phi (a_j,\la)\cos (\la(x-a_j))+\frac{\phi '(a_j^+,\la)}{\la }\sin (\la (x-a_j)), \ x\in (a_j,a_{j+1}).\label{soln}
\end{equation}
The details are left to the reader.
%\footnote{ For example,
%\begin{eqnarray*}
%\phi(a_{j+1},\la)& =& \phi (a_j,\la)\cos (\la \ell_j)+\frac{\phi '(a_j^+,\la)}{\la }\sin (\la \ell_j)\\
%& =& \sum_{n=0}^{j-1}b_n^j(\la)\la^n\cos (\la \ell_j)+\sum_{n=0}^{j}\tilde{b}_n^j(\la)\la^n\sin (\la \ell_j)\\
%& = & \sin (\la \ell_j)\tilde{b}_j^j\la^j+\sum_{n=0}^{j-1}(b_n^j(\la)\cos (\la \ell_j)+\tilde{b}_n^j(\la)\la^n\sin (\la \ell_j))\la^n.
%\end{eqnarray*}
%By the inductive hypothesis, there are at least $(n+1)$  sine terms in $(b_n^j(\la)\cos (\la \ell_j)+\tilde{b}_n^j(\la)\sin (\la \ell_j))$
%for each $n$.  }

We will now prove iii), again by induction. For $j=1$, the results hold by \rq{j=1}. Now assume the result
holds for $j<N+1$. By the inductive hypothesis and \rq{soln}, followed by a trigonometric identity, we have
\begin{eqnarray*}
\phi(a_{j+1},\la) & = & (\sin (\la a_j) +O(\la^3))\cos (\la \ell_j)+(\cos (\la a_j)+O(\la))\sin (\la \ell_j)\\
& = & \sin (\la a_{j+1})+O(\la^3).
\end{eqnarray*}
Similarly, we have
\begin{eqnarray*}
\phi'(a_{j+1}^+,\la)& = & \phi'(a_{j+1}^-,\la)-\la^2M_{j+1}\phi(a_{j+1},\la)\\
& = & \phi(a_j,\la)(-\la \sin (\la \ell_j))+\phi'(a_j^+,\la)( \cos (\la \ell_j))+O(\la^3)\\
& = & \la [(-\sin (\la a_j)\sin (\la (a_{j+1}-a_j)))+ \cos (\la a_j)\cos (\la (a_{j+1}-a_j))]+O(\la^3)\\
& = & \la \cos (\la a_{j+1})+O(\la^3).
\end{eqnarray*}
This completes the proof of iii).

The proof of iv is by induction. The details are left to the reader.

 \

In what follows, in the case $q=0$, we will denote the eigenvalue set as $\{ (\gamma_n)^2: n\in \bN \}$, while $\{ (\la_n)^2: n\in \bN \}$
will denote eigenvalues for $q\neq 0$. We define $\bK =\{ \pm 1, \pm 2, ...\}$. When $q=0$, the quadratic form $Q$ is positive definite, so
the eigenvalues $(\ga_n)^2$ are strictly positive. We choose $\ga_n>0$ for $n=1,2,...$, and then we can use
 the formula $\ga_{-n}=-\ga_n$, to extend the mapping $n \, \ra \ga_n$ from $\bN$ to $\bK$. We will refer to
$$\Gamma := \{ \ga_n:\ n\in \bK\}$$
as the eigenfrequencies of System \rq{Nmassq} for $q=0$.
We  define
$$G(\la )=\phi( \ell,\la).$$
Of course, $G(\la) =0$ whenever $\la =\ga_n$, and $\phi $ is a   corresponding eigenfunction.
We wish to examine some other  properties of $G$.
The following follows from Lemma \ref{q0asy}.
\begin{cor}\label{G}

\

a) The function $G(\la)/\la$ is an entire function of exponential type $\ell$ on the upper and lower half planes.

b) The roots of $G(\la)/\la$ are precisely $\Gamma$, including multiplicity.

c)  Let $\la =x+iy$, $x,y\in \bR$. Then there exists $y_0>0$ such that if $|y|>y_0$, then
there exist constants $C_0,C_1>0$ such that
$$C_0< \frac{|G(\la )|}{1+|\la|^N}< C_1, \  x \in \bR.$$

\end{cor}
Proof:

We begin with the proof of part b). It is clear by construction that $\la$ is a non-zero root of $G(\la )$ if and only if $\la =\ga_n$  for some $n\in \bK$.
 Regarding $\la =0$, we have $\ga_n \neq 0$ for any $n$,  but $G(0)=0$, and
by Lemma \ref{q0asy} part iii), $\lim_{\la \to 0}G(\la) /\la =\ell \neq 0$. Thus $\la$ is a root of $G(\la )/\la$ if and only if $\la =\ga_n$ for
some $n$.
Furthermore, the discussion of $G(\la )$ found in
 (\cite{Ince}, 10.72) shows that the zeros of $G(\la )/\la$ are all simple. In that work, where the equation $G(\la )/\la=0$ is called the
 ``characteristic equation", there are no masses, but the argument can easily be adapted
 to our setting.
 Part b) of the corollary now follows.
That $G(\la)/\la$ is an entire function also follows  from Lemma \ref{q0asy}.

To complete the proof of part a), let $\la =re^{i\theta}$. Fix $\theta$ and suppose $r$ is large. For $T(z)$ equalling either $\sin (z)$ or $\cos (z)$,
$$|T( \la \ell_j)|\sim e^{r\ell_j|\sin (\theta)|}.$$ Thus, using the notation of  Lemma \ref{q0asy} ,
$$\prod_{k=0}^N|T_k(\la \ell_k)|\sim e^{r\ell |\sin (\theta)|}.$$
Part a) now follows from \rq{bn} and \rq{Tk}.

To prove part c), we  use    Lemma \ref{q0asy} part i).
It follows that
$$\left|\frac{G(\la)}{\la^N} -(-1)^{N}\prod_{k=1}^{N}M_k \prod_{k=0}^{N}\sin (\la \ell_{k} )\right|=O(\la^{-1}),\  |\la | \gg 1.$$
Next, it is easy to show that there exists $y_0>0$ such that if $|y|>y_0$,
$$\prod_{k=0}^{N}|\sin (\ell_{k}(x+iy) )|\asymp 1,$$
this estimate being uniform in $x \in \bR$. The desired inequalities now follow easily.

\

 The following result gives, among other things, the asymptotics of $\Gamma $.
 \begin{cor}\label{spectrum}
Let $\Gamma '$ be any subset of $\Gamma$ obtained by deleting $2N$ elements. Then $\Gamma '$ can be reparametrized as
$$\Gamma ' = \cup_{j=0}^N\{ \ga^{(j)}_m\}_{m\in \bK},$$
where for each $j$,
 \begin{equation}
 |\ga_m^{(j)}-\frac{\pi m}{\ell_j}|=O(|m|^{-1}).\label{babyasy}
 \end{equation}
 \end{cor}

 {\bf Remarks}

 1- The case $N=M=1$, $q=0$ and $a=\ell/2$, studied in \cite{HZ}, shows that this result  is sharp
 in the sense that the exponent in the remainder term $O(|m|^{-1})$ cannot be improved in the general case.

 \

 2- The set $\cup_{j=0}^N\{ \frac{\pi m}{\ell_j}:\  m\in \bK \}$ is the frequency spectrum of
 the direct sum of Dirichlet  Laplacians on
 the intervals $(a_j,a_{j+1})$, with $j$ running from $0$ to $N$, or equivalently, the frequency spectrum  for  \rq{Nmassq} when $q=0$ and $M_j=\infty $
for each $j$.

 \

\noindent{\bf Proof of Corollary 2:}
Using Lemma \ref{q0asy} part i) and $u(\ell)=0$,  we get that
${\Gamma}$ is the solution set of the equation
 $$\la^{N-1}\prod_{k=0}^{N}\sin (\la \ell_k)=(-1)^{N+1}\frac{1}{\prod_{k=1}^NM_k}\sum_{n=0}^{N-1}b_n^{N+1}(\la){\la^{n-1}}.$$

We recall Rouch\'{e}'s Theorem.  Let $\Omega$ be a   planar domain bounded by a closed, simple curve $\partial \Omega$. Suppose $f(\la )$ and $g(\la)$ are analytic functions in an open set containing the closure of $\Omega$.  Suppose $|g(\la )|<|f(\la )|$ on $\partial \Omega$. Then $f$ and $f+g$ have the same number of zeros, including multiplicities, in $\Omega$.

We apply Rouch\'{e}'s Theorem as follows. Let
$$f(\la) = \la^{N-1}\prod_{k=0}^{N}\sin (\la \ell_k),\mbox{ and } g(\la )=(-1)^{N}\frac{1}{\prod_{k=1}^NM_k}\sum_{n=0}^{N-1}b_n^{N+1}(\la)
{\la^{n-1}}.$$
Thus ${\Gamma}$ is the zero set of $f+g$.
The zeros of $f$ are
$$\tilde{\Gamma}:
= \left ( \cup_{j=0}^N\{ \frac{\pi m}{\ell_j} : \ m\in \bK \}\right )\cup \{ 0\}  ,$$
with $0$ being of multiplicity $2N$.
We parametrize the zeros as
 $\{ \tilde{\gamma}_p: \ p\in \bZ\}$, with $\tilde{\gamma}_p$ listed in non-decreasing order and $\tilde{\gamma}_0=0$.
For $i\in \bK$, let $D_i$ be the disk of radius $C/|i|$ centered at $\tilde{\gamma}_{i}$; here $C$ is a positive, possibly large, constant to be given below.
 Let $\{ \Omega_m: \ m\neq 0 \} $ be the connected components of $\cup D_i$. We can choose the parametrization $m\ra \Omega_m $ so
 that $\Omega_{-m}=\{ \la :\  -\la \in \Omega_m\}$ and
 $$m_2>m_1 \mbox{ if and only if } \Re (z_2)>\Re (z _1) ,\ \forall (z_1,z_2)\in \Omega_{m_1}\times \Omega_{m_2}.$$

\underline{Part i: frequencies near $\pm \infty$.}

We will show that for $|p|>>1$, the  $\ga_p$ are all contained in $\Omega_m$ for some $m$.
 There exists $P_1=P_1(C)$ such that $|m|>P_1$ implies $\# \{ \tilde{\ga}_p \in \Omega_m\}\leq N+1.$ In what follows in part i, we will assume $|m|>P_1$.
 As a consequence,  there exist constants $C_i$,  independent of $(m,\la ,p)$ such that
 $$\la \in \Omega_{|m|} \mbox{ if and only if } C_0|m|\geq |\la |\geq C_1|m|,\mbox{ and } \tilde{\ga }_{|p|} \in \Omega_{|m|} \mbox{ if and only if }C_2|m|\leq |p|\leq C_3|m|.$$

We  estimate $f,g$ on  $\partial \Omega_m$ for some fixed $m$. We assume for now $p,m>0$; the proof for the other case is the same.
Throughout this paragraph, we assume $\la$ is in the closure of $ \Omega_m$.
 Denote  the zeros of $f$ within $\Omega_m$ as  $\tilde{\ga}_{p_m},..., \tilde{\ga}_{p_m+r}$.
The distance from $\partial \Omega_m$ to the set $\{ \tilde{\gamma}_p\}$ is at least $C/p_m$.
 Since $r\leq N+1$,  the diameter of $\Omega_m$ is less
than $2(N+1)C/p_m$.
Thus for each $k$ such that $\sin (\la \ell_k)$ vanishes in $\Omega_m$,
we have by Taylor's Theorem that
 $$\frac{C\ell_k}{2p_{m}}<|\sin (\la \ell_k )|<\frac{2(N+1)C\ell_k}{p_{m}}.\ \la \in \partial \Omega_m$$
For each $k$ such that $\sin (\la \ell_k)$  does not vanishes in $\Omega_m$,
by the boundedness of the derivative of the sine function we have
we have
 $$\frac{L_k}{p_{m}}< |\sin (\la \ell_k )|<\frac{2(N+1)L_k}{p_{m}},\ \la \in \partial \Omega_m$$
 where $L_k>\frac{C\ell_k}{2}.$ We remark that $L_k$ could be large and will depend on $m$, but that won't matter in
 the argument that follows.
   In what follows, it will be convenient to set $r=N+1-l$, with $0\leq l\leq N$. For simplicity of presentation, we will assume
   the $k$ for which $\sin (\la \ell_k)$  does not vanishes in $\Omega_m$ are $k=0,...,l-1$. The argument can easily be adapted to the other cases.
   In what follows, we will assume without loss of generality that $C>1$ and $L_k>1$ for $k=0,...,l-1$.
    Also, we will denote various constants  independent of $(m,\la, p, C)$
  by $\kappa_j$.
We have
$$
|f(\la )|\geq  \kappa_0 \left ( \frac{C}{p_{m}}\right )^{N+1-l}\left ( \prod_{k=0}^{l-1}\frac{L_k}{p_{m}}\right )
 |\la |^{N-1}\geq \kappa_1 \left ( \frac{C}{|\la|}\right )^{N+1-l}\left ( \prod_0^{l-1}\frac{L_k}{|\la |}\right )|\la |^{N-1}
 $$
\begin{equation}
 = \kappa_1 C^{N+1-l}\left ( \prod_0^{l-1}{L_k}\right ) |\la|^{-2}.\label{rou1}
\end{equation}
By Lemma \ref{q0asy} part i, $b_n^{N+1}$ can be written as a finite sum of terms, each which is a constant times $\prod_{k=0}^NT(\la \ell_k)$,
where  $s$ of the $T(z)$ are $\sin (z)$, with
$$ s\geq n+1.$$
%Let $Z$ be the number of such terms for $n$ running from 0 to $N-1$.
 We wish to show that for $C$  sufficiently large,
$\kappa_1 C^{N+1-l}(\prod_0^{l-1}{L_k})\ |\la|^{-2}$ is much larger than each of these terms, so that
\begin{equation}
\kappa_1 C^{N+1-l}\left ( \prod_0^{l-1}{L_k}\right ) |\la|^{-2}>\frac{1}{\prod_{k=1}^NM_k}\sum_{n=0}^{N-1}|b_n^{N+1}(\la)|{|\la |^{n-1}}
=|g(\la)|.\label{sticky}
\end{equation}
Consider one such term, which we label $\al \prod_{k=0}^N\hat{T}(\la \ell_k)$, with $\al $ a constant.
% Thus for each $k$ such that $\sin (\la \ell_k)$ vanishes on $\Omega_m$,
%we have
%by Taylor's Theorem, $|\sin (\la \ell_k)|< \frac{4(N+1)C\ell_k}{p_m}$ for $\la \in \partial \Omega_m$. For all other $k$ we have
% $$|\sin (\la \ell_k )|>\frac{4L_k}{p_{m}}.$$

There are four possible cases:

Case i:  $s=n+1$, $l<n+1$.

Since $n+1\leq N$, we have $s\leq N$, so
for $\la \in \partial \Omega_m$ we have
$$
|\al \prod_{k=0}^N\hat{T}(\la \ell_k) {\la^{n-1}}|\leq
\kappa_3  \big (\frac{C}{p_m} \big )^{s-l}\left ( \prod_{k=0}^{l-1}\frac{L_k}{p_{m}}\right ){|\la|^{n-1}}$$
\begin{equation}\lab{rou2}
\leq  \kappa_4 \big (\frac{C}{|\la |} \big )^{s-l}\left ( \prod_{k=0}^{l-1}\frac{L_k}{|\la |}\right ){|\la|^{n-1}}
\leq   \kappa_4  {C}^{N-l}\left ( \prod_{k=0}^{l-1}{L_k}\right )|\la|^{-2}.
\end{equation}
%Suppose first that $s=n+1$, so
% $$|\al \prod_{k=0}^N\hat{T}(\la \ell_k) \la^{n-1}|\leq  \kappa_2  {C}^{n+1}|\la|^{-2}.$$

Case ii: $s=n+1$, $l\geq n+1$.

In this case, $ \prod_{k=0}^N\hat{T}(\la \ell_k)$ might not vanish in $\Omega_m$, but
we can argue as follows:
\begin{equation}
|\al \prod_{k=0}^N\hat{T}(\la \ell_k) {\la^{n-1}}|\leq
  \kappa_5 \left ( \prod_{k=0}^{s-1}\frac{L_k}{|\la |}\right ){|\la|^{n-1}}
\leq \kappa_5\left ( \prod_{k=0}^{l-1}{L_k}\right ) |\la |^{-2}.\label{rou3}
\end{equation}

Case iii:  $N+1\geq s>n+1$,  $l<n+1$.

In this case, the argument in Case i yields
\begin{equation}
|\al \prod_{k=0}^N\hat{T}(\la \ell_k) {\la^{n-1}}|\leq \kappa_4  {C}^{N-l}\left ( \prod_{k=0}^{l-1}{L_k}\right )|\la|^{-3}.\label{rou4}
\end{equation}

Case iv:  $N+1\geq s>n+1$,  $l\geq n+1$.

In this case, the argument in Case ii yields
\begin{equation}
|\al \prod_{k=0}^N\hat{T}(\la \ell_k) {\la^{n-1}}|\leq \kappa_5 \left ( \prod_{k=0}^{l-1}{L_k}\right )|\la|^{-3}.\label{rou5}
\end{equation}

We now show that the terms in $g$ associated with Cases i, ii are dominated by the left hand side of \rq{sticky}.
Comparing \rq{rou1} with \rq{rou2} and \rq{rou3}, it is clear that we can choose $C$ (and increase $P_1$ if necessary) so that the sum of these terms in $|g|$  is  smaller than
$|f(\la)|/3$.
We now fix such $C$. For the terms in Cases iii, iv,   we see by comparing \rq{rou4} and \rq{rou5} with \rq{rou1} that
 by increasing $P_1$  if necessary,
the sum of such terms will be smaller than  $\frac{1}{3}|f(\la )|$.
Thus   $|f|>|g|$ on $\partial \Omega_m$, and so  Rouch\'{e}'s Theorem applies.

\underline{Part ii: frequencies near 0.}

We now locate, and implicitly count, the frequencies in $\tilde{\Gamma}$ that are ``near'' 0.
Let $L_1=\inf \{ |\la |: \ \la \in  \cup_{|m|\geq P_1}\Omega_m\} .$
Let $\tilde{\Omega}_L=\{ \la : |\Re (\la )|\leq L, |\Im \la |\leq L\}  $, with $L>L_1$ to be given below.
We  estimate $f,g$ on $\partial \tilde{\Omega}_L$.
On one hand, there exists $\ep >0$ such that
one can choose arbitrarily large  $L$ so that the lines $|\Re \la| =L$ do not intersect $\cup_{|m|\geq P_1} \Omega_m$, and furthermore
for $|\Re \la| =L$ we have
$$|\sin (\la \ell_k )|>\ep {e^{\ell_k |\Im \la |}}, \ k=0,..., N,$$
and hence
$$|f(\la )|> \ep^{N+1}{e^{\ell |\Im \la |}}L^{N-1}, \ \forall \la \in \partial \tilde{\Omega}_L.$$
On the other hand,
 by Lemma \ref{q0asy} part i, there exists a constant $C_5>0$ such that for each $n$ we have
$$|b^{N+1}_n(\la )|< C_5 e^{\ell |\Im \la |}, \mbox{ for } \la \in \partial \tilde{\Omega}_L, $$
and hence there exists $C_6>0$ such that
$$|g(\la )|=\frac{1}{\prod_{k=1}^NM_k}\sum_{n=0}^{N-1}|b_n^{N+1}(\la)|{|\la |^{n-1}}< C_6 {e^{\ell |\Im \la |}}L^{N-2}, \mbox{ for } \la \in \partial \tilde{\Omega}_L. $$
In particular, we can choose $L>L_1$ sufficiently large that $|f(\la )|>|g(\la )|$ on $\partial \Omega_0:= \partial \tilde{\Omega}_L$, so
Rouch\'{e}'s Theorem applies there.
Let $P_2=\min \{ m: |\la| >L,\  \forall \la \in \Omega_m\}.$ Then another argument using Rouch\'{e}'s Theorem
can be used to prove $all$ of $\Gamma$ are contained in $\Omega_0 \cup \big (\cup_{P_2}^{\infty}\Omega_m\big )$;
the details of this are left to the reader.
Since  $ \tilde{\Gamma}$ is contained in this set, we have established a bijection between the sets
${ \Gamma}$ and $ \tilde{\Gamma}$,
with the asymptotics announced \rq{babyasy}. The corollary now follows.

\

\

We now consider the asymptotics of system \rq{Nmassq} with $q\neq 0$.
For this paragraph, we paramatrize the sets
$\{ (\la_n)^2:\ n=1,2,...\}$ and $\{ (\ga_n)^2:\ n=1,2,...\}$, both in increasing order.
We  now define the eigenfrequencies $\{ \la_n: \ n\in \bK \}$ associated to $q\neq 0$; all of these appears by pairs.
If the $(\la_n)^2>0$, the associated pairs of  frequencies are defined by $\la_n$ and $\la_{-n}$, where for $n>0$ we set
$\la_n=|\la_n|$, and $\la_{-n}=-\la_n$.
For $(\la_n)^2<0$ the associated pair will be the complex conjugate pair $\la_{\pm n}=\pm i\sqrt{|\la_n|^2}$, and
for $\la_n^2=0$, we have $0=\la_n=\la_{-n}$.
We now compare eigenfrequencies
$\La = \{ \la_n: \ n\in \bK \}$ associated to $q\neq 0$ with those associate to $q=0$, namely  $\Gamma=\{ \gamma_n: n\in \bK \}$.
\begin{prop}\label{qasymp}
There exists a constant $C>0$ such that
$$ |\la_n-\ga_n|< \frac{C}{|n|},\ \forall n\in \bK.$$
\end{prop}
Proof: The quadratic forms associated to the system \rq{Nmassq} with $q=0$ and $q \neq 0$ are
$$Q_0(u)=\sum_{k=0}^N\int_{x=a_k}^{a_{k+1}}|u'(x)|^2+\sum_{k=1}^NM_k|u(a_k)|^2,$$
and $Q(u)=Q_0(u)+\int_0^{\ell}q(x)|u(x)|^2dx,$
with  quadratic form domain  for both $Q$ and $Q_0$ being ${\cal Q}$.
Suppose $K_0,K_1$ satisfy
$K_0<q(x)<K_1$ for  all $x$. Then for all $u$ in ${\cal Q}$,
$$Q_0(u)+K_0\leq Q(u)\leq Q_0(u)+K_1.$$
Using a standard mini-max argument (see, e.g. \cite{LL}), we get
$$\ga_n^2+K_0\leq \la_n^2\leq \ga_n^2+K_1,\ \forall n.$$
The proposition easily follows.

\

We conclude this section by remarking that Theorem \ref{spectrumq} follows immediately from Corollaries \ref{spectrum} and \ref{qasymp}.

\end{subsection}
\begin{subsection}{Riesz sequences}
In Section 5, we will prove Theorem \ref{mainthm} by solving a certain moment problem. This approach is complicated by the fact that, by Corollary \ref{spectrum} and Proposition \ref{qasymp}, the frequencies $\{ \la_n\}$
  have no uniform spectral gap. As a consequence, the exponential family $\{ e^{i\la_n t}\}$ is not a {\it Riesz sequence} (i.e. a Riesz basis in the closure of its linear span)  in $L^2(0,T)$ for any $T>0$. The purpose of this subsection is to construct the needed Riesz sequence using exponential divided differences.

 {\bf Definition}  A countable set of complex numbers $\{\mu_n\}$ called {\it uniformly discrete}
if \\  $\inf_{n \neq m}|\mu_n - \mu_m| >0.$ A set is {\it relatively uniformly discrete} if it is a union of a finite number of uniformly discrete sets.

{\bf Definition}
Assume $\{ \mu_j\}$ is a non-repeating sequence.
The exponential divided difference (E.D.D.) of order zero for $\{ e^{i\mu_n t}\}$ is $[e^{i\mu_1 t} ](t):=e^{i\mu_1t}.$
The E.D.D. of order $n-1$ is given by
$$[e^{i\mu_1 t},\ldots , e^{i\mu_n t}]=\frac{[e^{i\mu_1 t},\ldots ,e^{i\mu_{n-1} t}]-[e^{i\mu_2 t},\ldots ,e^{i\mu_n t}]}{\mu_1-\mu_n}.$$

One then easily derives the formula
$$
[e^{i\mu_1 t},\ldots ,e^{i\mu_n t}]=\sum_{k=1}^n\frac{e^{i\mu_kt}}{\prod_{j\neq k}(\mu_k-\mu_j)}.
$$
It is shown in \cite{AI1} that the functions $[e^{i\mu_1 t}]$, ..., $[e^{i\mu_1 t},\ldots ,e^{i\mu_n t}]$ depend on the parameters $\mu_j$ continuously, and each is invariant if we change the order of the $\mu_j$.

We will create E.D.D. from the set $\Lambda$, defined in the previous subsection.
We cite the following facts from \cite{AI1}.
 For any $z\in \bC$, denote by $D_z(r)$ the disk with center $z$ and radius $r$. Let $G^{p}(r)$, $p\in \bK$ be the connected components
 of the union $\cup_{z\in \Lambda}D_z(r)$. Write  $\Lambda^{p}(r)$ for the subset of $\Lambda$ lying in $G^{p}$, $\Lambda^{p}:
=\{ \la_n|\la_n\in G^{p}(r)\}$.
By Corollary \ref{spectrum}, $\Lambda$ can be decomposed into a  union of $N+1$ uniformly discrete sets, which we label $\Lambda_{j}$
with $j=0,..., N$.  Let
$$\delta_j:=\inf_{\la \neq \mu ;\la,\mu \in \Lambda_{j}}|\la -\mu|,\ \ \delta:=\min_j\delta_j.$$
Then for $r<r_0:=\frac{\delta}{2N+2}$, the number ${\cal N}^{(p)}(r)$ of elements of $\Lambda^{p}$ is
at most $N+1$.

It is now convenient to reparametrize $\Lambda$ according to the clusters $\{ \Lambda^{p}\}$:
$$\Lambda=\{ \la_l^p\ :\ p\in \bK, \ l=1,...,{\cal N}^{(p)}, \ \la_l^p\in  G^{p}(r)\}.$$
%\footnote {S: [ Confusing notations: first, several lines above $p$ run over $\mathbb{N}.$ Second,
%	it is better when $k\in \bK.$ Third, we should try to make notations similar to $ \gamma^{(j)}_n $ to help the reader.]
%J: [I fixed the first point. For the second point, I changed $k$ to $m$. I have not yet changed the third point, and
%I am not totally sure what you are saying here. I could alter $ \gamma^{j}_n $ to, say, $ \gamma_{j,n}$ if you think that
%helps, but such a  change will cause many changes in the next section. Do you have any suggestions?]}

We now fix $r<r_0$, and hence $G^p$ is independent of $r$ for each $p$.

%{\bf Remark on notation:} the reader should distinguish between $\{ \la^{(j)}_m\} $ from the previous subsection, where the frequencies have been %partitioned by its association with the subinterval
%$(a_j, a_{j+1})$, and this subsection's $\{ \la^p_l\} $, where the frequencies have been partitioned according to the clustering.

Denote by $\E$ be the family of  E.D.D. associated to $\Lambda$:
$$\E:=\cup_{p\in \bK}\cup_{k=1}^{{\cal N}^p}\{ [e^{i\la_1^pt}],[e^{i\la_1^pt},e^{i\la_2^pt}],\ldots , [e^{i\la_1^pt},\ldots ,
e^{i\la_{{\cal N}^{(p)}}^pt}] \}.$$
The following result follows from our Corollary \ref{spectrum}, Proposition \ref{qasymp}, and  Theorem 3 of \cite{AM1}, (also see \cite{AI1})
\begin{prop}\label{moran}
For any $T>2\ell$, $\E$ forms a Riesz sequence in $L^2(0,T)$.
\end{prop}

\end{subsection}
\begin{subsection}{More on spectral asymptotics}

Let $\{ \la_n^2 \}$  be the  eigenvalues for  the system \rq{Nmassq}, listed in increasing order,
and $\{ {\f}_n(x)\}_1^{\infty}$ the
 eigenfunctions normalized by
$||\f_n||_{L^2_M}=1$, where the $L^2_M$ norm is given by \rq{L2M}.
We conclude this section with a brief discussion on the asymptotics of ${\f}_n'(0)$. The term ${\f}_n'(0)$ appears in a moment problem associated with Theorem \ref{mainthm},
see \rq{bntN}, and hence the asymptotics  of $\f_n'(0)$ help to characterize the reachable set.

In the absence of any mass  it is well known that
\begin{equation}
|{\f}_n'(0)| \asymp |\la_n| + 1 \asymp n.\label{niceasy}
\end{equation}
These asymptotics also hold for the special
case of one mass when $a_1=\ell /2$ and $q=0$, which is the case where Hansen and Zuazua in \cite{HZ} used spectral methods to study controllability.
However,
in the typical case of attached masses, the first part of \rq{niceasy} will not hold.

 We will illustrate this here where
 we set $N=1,\, M=1$, $q=0$, $\ell =1$, and $a_1=a$ is arbitrary. For the purposes of our calculation,
  non-normalized eigenfunctions $\phi_n$ are chosen so that for $x\in (0,a)$ we have $\phi_n(x)=\sin (\la_n x)$, so that $\phi_n'(0)=\la_n.$
 Since
 $${\f}_n'(0)=\frac{\la_n}{\| \phi_n\|_{L^2_M}},$$
we will now compute $\| \phi_n\|_{L^2_M}$.
Easily we see
 $\int_0^a |\phi_n|^2\asymp a/2$ and $(\phi_n(a))^2=(\sin (\la_na))^2$.
  By \rq{j=1} and \rq{soln},
\begin{eqnarray}
\int_a^1 |\phi_n|^2 & = & \int_a^1|\sin (\la_n a)\cos(\la_n (x-a))+\big ( \cos (\la_n a)-\la_n\sin (\la_n a)\big )\sin (\la_n (x-a))|^2dx\nonumber \\
&=&\int_a^{1}|\sin (\la_n x)- \la_n \sin(\la_n a)\sin(\la_n(x-a))|^2dx.\label{efasy}
\end{eqnarray}
By Corollary \ref{spectrum} and Proposition \ref{qasymp}, the frequency spectrum with $2N$ terms deleted can be decomposed as  $\{ \la_n^{(0)}\}\cup \{ \la^{(1)}_n\}$, with
$$|\la_n^{(0)}- \frac{\pi n}{a}|=O(1/|n|)\mbox{ and }
|\la_n^{(1)}- \frac{\pi n}{1-a}|=O(1/|n|).$$
Denote the  corresponding eigenfunctions $\{ \phi_n^0 , \phi_n^1\}$.
If $\phi_n =  \phi_n^0$, then $|\la_n^{(0)}\sin (\la_n^{(0)}a)|=O(1)$, so
the integral in  \rq{efasy} is bounded as a function of $n$. Hence
$$|({\f}_n^0)'(0)|\asymp |\la_n^{(0)}|\asymp |n|,$$
same as \rq{niceasy}.

If $\phi_n=\phi^1_n$, we must consider two cases.

Case i: $\frac{a}{1-a}\in \bN$. This  condition is equivalent to  $a\in \{ 1/2, 2/3, 3/4, ...\}$. In these cases,
$|\sin (\la^{(1)}_n a)|=O(1/|n|)$.
Thus for  $\phi_n =  \phi_n^1$, then the integral in  \rq{efasy} is bounded as a function of $n$, and hence
$$|({\f}_n^1)'(0)|\asymp |\la_n^{(1)}|\asymp |n|,$$
same as \rq{niceasy}.

Case ii: Assume now $a \not \in \{ 1/2, 2/3, 3/4, ...\}$, or equivalently,
$\frac{a}{1-a}$ is not an integer. The proof of the  following result, an elementary exercise in analysis, is left to the reader.
\begin{lemma}
Suppose $\frac{a}{1-a}$ is not an integer. Let $\ep \in (0,\frac{1}{100})$. There exists  an infinite subsequence of $\bN$, labeled $\{ n_m: m=1, 2, ..,\}$ and listed in increasing order,
such that
$$|\la_{n_m}^{(1)}a-\pi p|>\ep,\ \forall p,m\in \bN ,$$
and such that   $n_m\asymp m$.
\end{lemma}
From the lemma, it follows that
 there exists constant $\delta >0$ such that
$$|\sin (\la^{(1)}_{n_m}a)|>\delta , \forall m\in \bN.$$
Thus by \rq{efasy}
 $$\int_a^1|\phi^1_{n_m}|^2 \asymp m,$$
and hence
$$|({\f}^1_{n_m})'(0)|\asymp 1.$$

\end{subsection}
\end{section}

\begin{section}{Proof of full controllability}

In this section we prove Theorem \ref{mainthm}.
 In what follows, we denote the frequencies of the system \rq{Nmassq} by $\La =\{ \la_k^p: p\in \bK, k=1,..., {\cal N}^{(p)}\}$,
 with   $\f_k^{|p|}$ the corresponding eigenfunctions,
 {orthonormal} in $L^2_M$. We will assume for simplicity  that the eigenvalues satisfy $(\la_k^p)^2>0$. If this were not the case in what follows, it would
 suffice to replace $\sin (\la t)/\la$ by  $\sinh (|\la|t)/|\la|$ in the case $\la^2<0$, and by $t$ in the case  $\la^2=0$.

\noindent{\bf Step 1: Fourier representation of solution to initial boundary value problem.}

 We present the solution of System \rq{A_N11}-\rq{control9}
in the form of the series
\beq \lab{uxtN}
u^f(x,t)=\sum_{p=1}^{\infty} \sum_{k=1}^{{\cal N}^{(p)}}a_k^p(t)\,\f_k^p(x).
\eeq

For any $T>0$ and $f\in L^2(0,T)$,
standard calculations using the weak solution formulation for System \rq{A_N11}-\rq{control9} (see, eg. \cite{AI} Ch. 3) give for each $p$,
\beq \lab{bntN}
a_k^p(t)=\frac{(\f_k^p)'(0)}{\la_k^p}\int_0^t f(\tau)\, \sin \la_k^p(t-\tau)\,d\tau, \, (a_k^p)'(t)={(\f_k^p)'(0)}\int_0^t f(\tau)\, \cos \la_k^p(t-\tau)\,d\tau .
\eeq

Denote $s_k^p(x)= \sin (\la_k^p x)$ and $c_k^p(x)=\cos (\la_k^px)$.
We set
$a_k^p=a_k^p(T),b_k^p=(a_k^p)'(T)$, and
\beq \lab{alphabeta}
\al_k^p=\frac{a_k^p\la_k^p}{(\f_k^p)'(0)}, \ \beta_k^p=\frac{b_k^p}{(\f_k^p)'(0)}.
\eeq
Let $<*,*>_T$ be the standard complex inner product on $L^2(0,T)$. Let $f^T(t)=f(T-t)$, so
 $\int_0^Tf(t)\sin ({\la (T-t)})dt=
<f^T,\sin (\la t)>_T$.
Then \rq{bntN} for $t=T$ can be written as
\beq \lab{an2}  \al_k^p=<f^T,s_k^p>_T, \ p  \in \bN,\ k=1, ..., {\cal N}^{(p)} \eeq
and
\beq \lab{bnN} \beta_k^p=<f^T,c_k^p>_T, \ p  \in \bN,\ k=1, ..., {\cal N}^{(p)}. \eeq

Using $e^{it}=\cos (t)+i\sin (t)$, we get the following equations which hold for all $ p  \in \bN,\ k=1, ..., {\cal N}^{(p)}$:
\begin{eqnarray*}
i\al_k^p+ \beta_k^p & = & <f^T,e^{-i\la_k^pt}>_T,\\
-i\al_k^p+ \beta_k^p & = & <f^T,e^{i\la_k^pt}>_T.
\end{eqnarray*}
We recall   $\la_{k}^{-p}=-\la_k^p$ for $p\geq 1$. Similarly
we set $\al_k^{-p}=-\al_k^p$, and
 $\beta_k^p=\beta_{k}^{|p|}$ for all $p\in \bK$.
Define $\ga_k^p$ by
\beq\lab{gamma_k}
\ga_k^p=
(- i\al_k^p+ \beta_k^p) ,\ \forall\ p  \in \bK,\ k=1, ..., {\cal N}^{(p)}.
\eeq
Then we have
\begin{equation}
<f^{T},e^{i\la_k^pt}>_{T}=\ga_k^p,\ p  \in \bK,\ k=1, ..., {\cal N}^{(p)}.\label{momentN}
\end{equation}
Assigning terminal data to System \rq{A_N11}-\rq{control9} is equivalent to assigning values to $\{ \gamma_k^p\}$, in which case
\rq{momentN} can be viewed as a moment problem. However, solving this moment problem is complicated by the fact that
$\{ e^{i\la_k^pt}\}$ does not form a Riesz sequence. Thus we shall use E.D.D. to
rewrite \rq{momentN}.

\

\noindent{\bf Step 2: some Riesz sequences.}

Recall $$\E=\cup_{p\in \bK}\{ [e^{i\la_1^pt}],[e^{i\la_1^pt},e^{i\la_2^pt}],\ldots , [e^{i\la_1^pt},\ldots ,e^{i\la_{{\cal N}^{(p)}}^pt}] \}.$$
We define
${\cal S}$ and ${\cal C}$ to be   corresponding families of divided differences of $s_n$ and $c_n$. Thus,
$${\cal S}=\cup_{p\in \bN} \{ [s_1^p(t)],\ldots , [s_1^p(t),\ldots ,s^p_{{\cal N}^{(p)}}(t)]\},$$
$${\cal C}=\cup_{p\in \bN} \{ [c_1^p(t)],\ldots , [c_1^p(t),\ldots ,c^p_{{\cal N}^{(p)}}(t)]\}.$$
%Recall that if $\sigma_e$ is a $N-1$ element subset of $\E$, then $\E  \setminus \sigma_e$
%forms a Riesz basis on $L^2(0,2\ell )$. It will be convenient to choose $\sigma$ as follows.

For any $N\geq 1$, we have that $\E$ is a Riesz sequence on $L^2(0,T)$ for $T>2\ell$
by Proposition \ref{moran}. It
 then follows from
  (\cite{ABI}, Lemma 5.1) that ${\cal S} $ and ${\cal C} $ form Riesz sequences
in $L^2(0,T)$ for $T>\ell$.

\

\noindent{\bf Step 3: a technical lemma.}

We wish to   rewrite the moment problem above in terms of our E.D.D., but first we need to develop some notation.
For $\la_k^p$ as above and $a_1^p, ..., a_n^p \in \bC$, we construct divided differences of these $numbers$ iteratively
by $[a_1^p]'=a_1^p$, and
$$[a_1^p,\ldots , a_n^p]'=\frac{[a_1^p,\ldots ,a_{n-1}^p]'-[a_2^p,\ldots ,a_n^p]'}{\la^p_1-\la^p_n}.$$
It is easy to see that
\begin{equation}
[a_1^p,\ldots ,a_n^p]'=\sum_{k=1}^n\frac{a_k^p}{\prod_{j\neq k}(\la^p_k-\la^p_j)}.\label{gdd}
\end{equation}

\begin{lemma}\label{twist}
For $n=1,...,{\cal N}^{(p)}$,
$$a^p_{n}=\sum_{k=1}^{n}[a_1^p,...,a_{k}^p]'\prod_{j=1}^{k-1}(\la^p_n-\la^p_j).$$
\end{lemma}
Here we use the convention that $\prod_{j=1}^{0}(\la^p_n-\la^p_j)=1$.

Proof:
For readability, we drop the $p$ superscript in what follows.
By \rq{gdd} we have
\begin{equation}
a_n=\prod_{j=1}^{n-1}(\la_n-\la_j)\left ( [a_1,...,a_n]'-\sum_{k=1}^{n-1}\frac{a_k/(\la_k-\la_n)}{\prod_{j\in \{ 1,...,(n-1)\} -k}(\la_k-\la_j)}\right ) .\label{anp}
\end{equation}
An algebra exercise shows that for any $m$,
$$\sum_{k=1}^{m}\frac{a_k/(\la_k-\la_n)}{\prod_{j\leq m;j \neq k}(\la_k-\la_j)}=\frac{1}{\la_m-\la_n}\left ( [a_1,...,a_{m}]'-
\sum_{k=1}^{m-1}\frac{a_k/(\la_k-\la_n)}{\prod_{j\leq (m-1);j \neq k}(\la_k-\la_j)}\right ) .$$
Applying this repeatedly to \rq{anp} gives the lemma.

\

\noindent{\bf Step 4: Solution to moment problem associated to full controllability}

It is easy to see that  \rq{momentN} is equivalent to
\beq \lab{momentEDD}
[\ga_1^p,\ldots , \ga_k^p]'=<f^{T},[e^{i\la_1^pt},\ldots , e^{i\la_k^pt}]>_{T};\ p\in \bK, \ k=1,..., {\cal N}^{(p)}.
\eeq
Let   $T_*>2\ell$. We now view \rq{momentEDD} as a moment problem. Let
%\footnote{this is where we need $T>2\ell$}
$[\ga_1^p,\ldots , \ga_k^p]' \in \ell^2$. Then
there exists a solution $f\in L^2(0,T_*)$ to \rq{momentEDD} with
\begin{equation}
\|  f\|^2_{L^2(0,T_*)}\asymp \sum_{p,k}|[\ga_1^p,\ldots , \ga_k^p]'|^2.\label{expR}
\end{equation}
Thus we have solved this moment problem, but it is not yet apparent how to characterize the reachable set of the functions corresponding to $\{ [\ga_1^p,\ldots , \ga_k^p]'\}  \in \ell^2 $. As a first step of this characterization,
 observe that for all $p\in \bK$, and all $k=1,..,{\cal N}^{(p)}$,
$$
[\beta_1^p,\ldots , \beta_k^p]' + i[\al_1^p,\ldots , \al_k^p]'=[\ga_1^p,\ldots , \ga_k^p]' .
$$
Hence
\beq
\sum_{p,k}|[\ga_1^p,\ldots , \ga_k^p]'|^2=\sum_{p,k}|[\al_1^p,\ldots , \al_k^p]'|^2+\sum_{p,k}|[\beta_1^p,\ldots , \beta_k^p]'|^2.\label{abq}
\eeq

%\noindent{\bf Step 2: We express the full control problem as a moment problem}

%We cannot solve this moment problem because it is not clear that $\{ \ga_j\} \in \ell^2$,
%{\bf S: [It is, certainly, clear for both: the Fourier coef. of the target data and for the Fourier coef. of the terminal state.]}
%{\bf J responds: I am not sure  you are right. Regarding the target data,
%if we assume only that $\{ y_{0,k}^p\} \in \ell^2$, then  the asymptotics of $\la_k^p$ and $(\f_k^p)'(0)$ discussed in Section 4 show that %$\al_k^p=\frac{y_{0,k}^p\la_k^p}{(\f_k^p)'(0)}$  need not be in $\ell^2$. I agree that the Fourier coeff. of the terminal state are in $\ell^2$, %because $\{ e^{i\la_k^pt}\}$ is a Bessel sequence  }
%and also $\{ e^{i\la_k^pt}\}$ is not a  Riesz sequence
%in $L^2(0,T)$
%for any given $T$. Thus we will restate \rq{momentN} in terms of EDD.

%where $\{ \ga_j\}$ is an appropriate space. To find the appropriate space, we have to study the properties
%of the exponential family $\{ e^{i\la_k^pt} \}$.

%In what follows, we will sometimes drop the superscript $p$ for readability, and for the same reason write $u$ for $u^f$. {\bf S: [ Not now. On Step % it will be fine, now it is better to keep the full notations.]}

\

\noindent{\bf Step 5: Moment problems associated to shape and velocity control}

We rewrite \rq{an2} in the form
\beq \lab{bnmN1}   \  [\al_{1}^p,...,\al_{k}^p]' =<f^T,[s_{1}^p,..., s_{k}^p]>_T,  \   \ p\in \bN, \ k=1,...,{\cal N}^{(p)}, \eeq
and \rq{bnN} in the form
\beq \lab{bnmN2}   \  [\beta_{1}^p,...,\beta_{k}^p]' =<f^T,[c_{1}^p,..., c_{k}^p]>_T,  \   \ p\in \bN, \ k=1,...,{\cal N}^{(p)}. \eeq
By Lemma \ref{twist}, dropping the superscript $p$ for readability,
\begin{eqnarray*}
a_1\f_1+\dots a_{\N}\f_{\N}& = & \sum_{n=1}^{\N}\f_n (\sum_{k=1}^{n}[a_1,...,a_{k}]'\prod_{j=1}^{k-1}(\la_n-\la_j))\\
& = & \sum_{k=1}^{\N}[a_1,...,a_{k}]'\ ( \sum_{j=k}^{\N}\f_j(\prod_{l=1}^{k-1}(\la_j-\la_l))).
\end{eqnarray*}
Thus,  with $\eta_k=\la_k/\f'_k(0)$ hence  $a^p_k=\al_k^p/\eta_k^p$,   we can rewrite \rq{uxtN}, with $t=T$, as
\begin{eqnarray}
u(x,T)& = & \sum_{p=1}^{\infty}\sum_{k=1}^{\N^{(p)}}a^{p}_k\f^{p}_k\nonumber\\
& = & \sum_{p=1}^{\infty}\sum_{k=1}^{\N^{(p)}}\al^{p}_k\f^{p}_k/\eta_k^p\nonumber\\
& = & \sum_{p=1}^{\infty}\sum_{k=1}^{\N^{(p)}}[\al_1^p,\ldots , \al_k^p]'(\sum_{j=k}^{\N^{(p)}}\frac{\f_j^p}{\eta_j^p}\prod_{l=1}^{k-1}(\la_j^p-\la_{l}^p)).\label{u0N}
\end{eqnarray}
%Let $\theta=\{ (k,p):\ p=1,..., Q; k=1, ..., \N^{(p)}-1\} \cup  \{ (k,p):\ p>Q; k=1, ..., \N^{(p)}\}.$
One can similarly express $u_t$ as a series in the terms $\sum_{j=k}^{\N^{(p)}}\frac{\f_j^p\la^p_j}{\eta_j^p}\prod_{l=1}^{k-1}(\la_j^p-\la_{l}^p)$.

\noindent{\bf Step 6:   Riesz bases for $W_0$ and $W_{-1}$.}

 Define
$$\psi^p_k: = \sum_{j=k}^{\N^{(p)}}\frac{\f_j^p}{\eta_j^p}\prod_{l=1}^{k-1}(\la_j^p-\la_{l}^p). $$
 We will  prove that the family $ \{\psi^p_k \ | \  p\geq 1, \ k=1,..., {\cal N}^{(p)} \} $
 forms a Riesz basis of $ W_0.$ This result is not only central to the proof of  Theorem 2, but is of independent interest.
Let  $v$ be an arbitrary element of $W_0$.  We need to prove the following facts:
 	
\noindent(i) $v$ can be presented uniquely as a  series with respect to $\psi_k^p$, with
 	quadratically summable coefficients and the series converging in $W_0$,

\noindent (ii) the $\ell^2-$norm of the coefficients is equivalent to $W_0-$norm of the function.

We begin proving part i.
 Suppose $T>\ell$.
  By Theorem \ref{shapecont}, part B,
there exists $g \in L^2(0,T)$ such that $u^g(x,T)=v(x)$,
with
\beq \lab{snorm}
||v||_{W_0}\asymp ||g||_{L^2(0,T)}.
\eeq
 We introduce simplifying notations and rewrite \rq{bnmN1} in the form
\beq \lab{aln*}
\nu_n = <g^T,S_n>_T, \   \ n  \in \bN,
\eeq
and rewrite \rq{u0N} in the form
\beq \lab{uxf*} v(x)= u^g(x,T)=
\sum_{n \in \bN} \nu_n \psi_n.
\eeq
Since $\f_n$ form an orthogonal basis in $L_M^2(0,\ell),$ it follows that $v $ can be
uniquely represented as series  \rq{uxtN}, and  hence as series  \rq{uxf*}.
Since  $\{ S_n:n\geq 1\}$ forms a Riesz sequence on $L^2(0,T)$,  it follows from \rq{aln*} and Theorem 3 that $\nu_n\in \ell^2$. Below we
 will prove the series \rq{uxf*} converges in $W_0.$

We now prove part ii above. Denote by ${\cal S}$ the closure of linear span of $\{S_n\}$ in $L^2(0,T)$. Let $g_0$ be the orthogonal projection of $g$ onto ${\cal S}$. Then clearly $g_0$
satisfies the moment problem \rq{aln*}, and hence $u^{g_0}(x,T)=v(x)$ and  by Proposition \ref{Nreg}
\beq \lab{precg}
||v||_{W_0}\prec ||g_0||_{L^2(0,T)}.
\eeq
On the other hand, by \rq{snorm} we have
$$||g_0||_{L^2(0,T)}\leq ||g||_{L^2(0,T)} \prec ||v||_{W_0}.$$
Thus
$$
||g_0||_{L^2(0,T)}\asymp ||v||_{W_0}.
$$
Furthermore, since $\{S_n\}$ form a Riesz basis on ${\cal S}$, we have
$$||g_0||_{L^2(0,T)}^2\asymp \sum_{n \in \bN} |\nu_n|^2.$$
Combining,  we get
$$
||\sum_{n \in \bN}\nu_n\psi_n||_{W_0}^2=||v||_{W_0}^2\asymp \sum_{n \in \bN} |\nu_n|^2.
$$
This argument also shows that for any sequence $\{ \eta_n: \ n\in \bN \}\in \ell^2$,
$$
||\sum_{n \in \bN}\eta_n\psi_n||_{W_0}^2\asymp \sum_{n \in \bN} |\eta_n|^2.
$$
This proves ii and that  series \rq{uxf*} converges in $W_0.$

We similarly find a Riesz basis for $W_{-1}$.
Set
$$\tilde{\psi}^p_k=\sum_{j=k}^{\N^p}{\f_j^p}(x)(\f_k^p)'(0)\prod_{l=1}^{k-1}(\la_j^p-\la_{l}^p)\ | \ p\geq 1, \ k=1,..., {\cal N}^{(p)}.$$
We have
$$u_t(x,T)= \sum_{p=1}^{\infty}\sum_{k=1}^{\N^{(p)}}[\beta_1^p,\ldots , \beta_k^p]'\tilde{\psi}_k^p(x), $$
and arguing as above, one can show
$$\{ \tilde{\psi}_k^p\ | \ p\geq 1, \ k=1,..., {\cal N}^{(p)} \}
\mbox{ forms a Riesz basis of }\ W_{-1}.$$
The details are left to the  reader.

\noindent{\bf Step 7: completion of proof of theorem}

Let $(y_0,y_1)\in W_0\times W_{-1}$.
We wish to find   $f\in L^2(0,T_*)$ solving
\begin{equation}
u^f(x,T_*)=y_0(x),\ u_t(x,T_*)=y_1(x).\label{control*}
\end{equation}
We begin by expressing \rq{control*} as a moment problem.
We expand $y_0$ and $y_1$ using the Riesz bases from Step 6:
$$y_0(x)=\sum_{p=1}^{\infty}\sum_{k=1}^{{\cal N}^{(p)}}y_{0,k}^p\psi_k^p(x),\ y_1(x)=\sum_{p=1}^{\infty}\sum_{k=1}^{{\cal N}^{(p)}}y_{1,k}^p\tilde{\psi}_k^p(x).$$
Extend $y_{0,k}^p$, resp. $y_{1,k}^p$ as even, resp. odd functions in $p\in \bK$.
Set $\hat{\ga}_k^p=y_{0,k}^p+iy_{1,k}^p\in \ell^2$.
 Then there exists $f\in L^2(0,T_*)$ solving  following moment problem:
\begin{equation}
<f^{T_*},[e^{i\la_1^pt},...,e^{i\la_{k}^pt}]>_{T_*}=\hat{\ga}_k^p,\ p  \in \bK,\ k=1, ..., {\cal N}^{(p)},\label{momentN1}
\end{equation}
and hence the same $f$ solves \rq{control*}.
We  have by Step 5 along with \rq{u0N} that
$$\sum_{p=1}^{\infty}\sum_{k=1}^{\N^{(p)}}|y_{0,k}^p |^2 \asymp \| y_0\|^2_{W_0}.$$
Similarly, Step 5 implies
$$\sum_{p=1}^{\infty}\sum_{k=1}^{\N^{(p)}}|y_{1,k}^p |^2 \asymp\| y_1\|^2_{W_{-1}}.$$
Thus by \rq{expR} and \rq{abq}, we get
$$\|  f\|^2_{L^2(0,T_*)}\asymp \sum_{p,k}|\hat{\ga}_k^p|^2
\asymp \| y_0\|^2_{W_0}+\| y_1\|^2_{W_{-1}}.$$
The proof is complete.

\end{section}

%\begin{section}{Concluding Remarks}

%\end{section}

%\begin{section}{References}

\end{document}